\documentclass[a4paper,10pt]{article}
\usepackage{amsmath,amsthm,amssymb,amscd,layout}   
\usepackage{bm,color,enumerate} 
\usepackage{amsfonts}
\usepackage{mathrsfs}
\usepackage[dvips]{graphicx}

\usepackage{tensor}
\usepackage{enumitem} 
\usepackage{graphicx,latexsym}
\usepackage{dsfont}
\usepackage{pgfplots}
\usepackage{tikz}
\usepackage{bm}
\usepackage{caption}

\newtheorem{theorem}{Theorem}[section]
\newtheorem{lemma}[theorem]{Lemma}
\newtheorem{proposition}[theorem]{Proposition}

\theoremstyle{definition}
\newtheorem{definition}[theorem]{Definition}
\newtheorem{example}[theorem]{Example}
\newtheorem{claim}[theorem]{Claim}

\theoremstyle{remark}
\newtheorem{remark}[theorem]{Remark}

\theoremstyle{conjecture}

\numberwithin{equation}{section}
\newcommand{\hilsp}{\mathscr{H}}
\newcommand{\drawarc}[2]{
\ifnum#1=1
  \ifnum#2=1
      \put(  0, 0){\line(0,1){5}}
  \fi
  \ifnum#2=2
      \qbezier(0, 0)( 4, 8)( 8, 0)
  \fi
  \ifnum#2=3
      \qbezier(0, 0)( 8,16)(16, 0)
  \fi
  \ifnum#2=4
      \qbezier(0, 0)(12,24)(24, 0)
  \fi
  \ifnum#2=5
      \qbezier(0, 0)(16,32)(32, 0)
  \fi
  \ifnum#2=6
      \qbezier(0, 0)(20,40)(40, 0)
  \fi
  \ifnum#2=7
      \qbezier(0, 0)(24,48)(48, 0)
  \fi
  \ifnum#2=8
      \qbezier(0, 0)(28,56)(56, 0)
  \fi
  \ifnum#2=9
      \qbezier(0, 0)(32,64)(64, 0)
  \fi
  \ifnum#2=10
      \qbezier(0, 0)(36,72)(72, 0)
  \fi
\fi
\ifnum#1=2
  \ifnum#2=2
      \put(  8, 0){\line(0,1){5}}
  \fi
  \ifnum#2=3
      \qbezier( 8, 0)(12, 8)(16, 0)
  \fi
  \ifnum#2=4
      \qbezier( 8, 0)(16,16)(24, 0)
  \fi
  \ifnum#2=5
      \qbezier( 8, 0)(20,24)(32, 0)
  \fi
  \ifnum#2=6
      \qbezier( 8, 0)(24,32)(40, 0)
  \fi
  \ifnum#2=7
      \qbezier( 8, 0)(28,40)(48, 0)
  \fi
  \ifnum#2=8
      \qbezier( 8, 0)(32,48)(56, 0)
  \fi
  \ifnum#2=9
      \qbezier( 8, 0)(36,64)(64, 0)
  \fi
  \ifnum#2=10
      \qbezier( 8, 0)(40,72)(72, 0)
  \fi
\fi
\ifnum#1=3
  \ifnum#2=3
      \put( 16, 0){\line(0,1){5}}
  \fi
  \ifnum#2=4
      \qbezier(16, 0)(20, 8)(24, 0)
  \fi
  \ifnum#2=5
      \qbezier(16, 0)(24,16)(32, 0)
  \fi
  \ifnum#2=6
      \qbezier(16, 0)(28,24)(40, 0)
  \fi
  \ifnum#2=7
      \qbezier(16, 0)(32,32)(48, 0)
  \fi
  \ifnum#2=8
      \qbezier(16, 0)(36,40)(56, 0)
  \fi
  \ifnum#2=9
      \qbezier(16, 0)(40,48)(64, 0)
  \fi
  \ifnum#2=10
      \qbezier(16, 0)(44,48)(72, 0)
  \fi
\fi
\ifnum#1=4
  \ifnum#2=4
      \put( 24, 0){\line(0,1){5}}
  \fi
  \ifnum#2=5
      \qbezier(24, 0)(28, 8)(32, 0)
  \fi
  \ifnum#2=6
      \qbezier(24, 0)(32,16)(40, 0)
  \fi
  \ifnum#2=7
      \qbezier(24, 0)(36,24)(48, 0)
  \fi
  \ifnum#2=8
      \qbezier(24, 0)(40,32)(56, 0)
  \fi
  \ifnum#2=9
      \qbezier(24, 0)(44,40)(64, 0)
  \fi
  \ifnum#2=10
      \qbezier(24, 0)(48,48)(72, 0)
  \fi
\fi
\ifnum#1=5
  \ifnum#2=5
      \put( 32, 0){\line(0,1){5}}
  \fi
  \ifnum#2=6
      \qbezier(32, 0)(36, 8)(40, 0)
  \fi
  \ifnum#2=7
      \qbezier(32, 0)(40,16)(48, 0)
  \fi
  \ifnum#2=8
      \qbezier(32, 0)(44,24)(56, 0)
  \fi
  \ifnum#2=9
      \qbezier(32, 0)(48,32)(64, 0)
  \fi
  \ifnum#2=10
      \qbezier(32, 0)(52,40)(72, 0)
  \fi
\fi
\ifnum#1=6
  \ifnum#2=6
      \put( 40, 0){\line(0,1){5}}
  \fi
  \ifnum#2=7
      \qbezier(40, 0)(44, 8)(48, 0)
  \fi
  \ifnum#2=8
      \qbezier(40, 0)(48,16)(56, 0)
  \fi
  \ifnum#2=9
      \qbezier(40, 0)(52,24)(64, 0)
  \fi
  \ifnum#2=10
      \qbezier(40, 0)(56,32)(72, 0)
  \fi
\fi
\ifnum#1=7
  \ifnum#2=7
      \put( 48, 0){\line(0,1){5}}
  \fi
  \ifnum#2=8
      \qbezier(48, 0)(52, 8)(56, 0)
  \fi
  \ifnum#2=9
      \qbezier(48, 0)(56,16)(64, 0)
  \fi
  \ifnum#2=10
      \qbezier(48, 0)(60,24)(72, 0)
  \fi
\fi
\ifnum#1=8
  \ifnum#2=8
      \put( 56, 0){\line(0,1){5}}
  \fi
  \ifnum#2=9
      \qbezier(56, 0)(60,40)(64, 0)
  \fi
  \ifnum#2=10
      \qbezier(56, 0)(64,48)(72, 0)
  \fi
\fi
\ifnum#1=9
  \ifnum#2=9
      \put( 64, 0){\line(0,1){5}}
  \fi
  \ifnum#2=10
      \qbezier(64, 0)(68,8)(64, 0)
  \fi
\fi
\ifnum#1=10
  \ifnum#2=10
      \put( 72, 0){\line(0,1){5}}
  \fi
\fi
}

\newcommand{\stos}{
   \put( 10, 0){\line(0,1){6}}
   \put( 10, 1){\circle*{2}}
}
\newcommand{\atoa}{
   \put( 0, 10){\line(1,0){20}}
}
\newcommand{\btob}{
   \put( 0, 18){\line(1,0){20}}
}
\newcommand{\ctoc}{
   \put( 0, 26){\line(1,0){20}}
}
\newcommand{\dtod}{
   \put( 0, 34){\line(1,0){20}}
}
\newcommand{\etoe}{
   \put( 0, 42){\line(1,0){20}}
}
\newcommand{\stoa}{
   \put( 19, 1){\oval(18,18)[tl]}
   \put( 19,10){\line(1,0){1}}
   \put( 10, 1){\circle*{2}}
}
\newcommand{\atob}{
   \qbezier( 0,10)( 5,10)(10,14)
   \qbezier(19,18)(15,18)(10,14)
   \put(19,18){\line(1,0){1}}
}
\newcommand{\btoc}{
   \qbezier( 0,18)( 5,18)(10,22)
   \qbezier(19,26)(15,26)(10,22)
   \put(19,26){\line(1,0){1}}
}
\newcommand{\ctod}{
   \qbezier( 0,26)( 5,26)(10,30)
   \qbezier(19,34)(15,34)(10,30)
   \put(19,34){\line(1,0){1}}
}
\newcommand{\dtoe}{
   \qbezier( 0,34)( 5,34)(10,38)
   \qbezier(19,42)(15,42)(10,38)
   \put(19,42){\line(1,0){1}}
}
\newcommand{\etod}{
   \qbezier( 0,42)( 5,42)(10,38)
   \qbezier(19,34)(15,34)(10,38)
   \put(19,34){\line(1,0){1}}
}
\newcommand{\dtoc}{
   \qbezier( 0,34)( 5,34)(10,30)
   \qbezier(19,26)(15,26)(10,30)
   \put(19,26){\line(1,0){1}}
}
\newcommand{\ctob}{
   \qbezier( 0,26)( 5,26)(10,22)
   \qbezier(19,18)(15,18)(10,22)
   \put(19,18){\line(1,0){1}}
}
\newcommand{\btoa}{
   \qbezier( 0,18)( 5,18)(10,14)
   \qbezier(19,10)(15,10)(10,14)
   \put(19,10){\line(1,0){1}}
}
\newcommand{\atos}{
   \put(  1, 1){\oval(18,18)[tr]}
   \put(  0,10){\line(1,0){1}}
   \put( 10, 1){\circle*{2}}
}
\newcommand{\btos}{
   \put(  1, 9){\oval(18,18)[tr]}
   \put(  0,18){\line(1,0){1}}
   \put( 10, 1){\line(0,1){8}}
   \put( 10, 1){\circle*{2}}
}
\newcommand{\ctos}{
   \put(  1,17){\oval(18,18)[tr]}
   \put(  0,26){\line(1,0){1}}
   \put( 10, 1){\line(0,1){16}}
   \put( 10, 1){\circle*{2}}
}
\newcommand{\dtos}{
   \put(  1,25){\oval(18,18)[tr]}
   \put(  0,34){\line(1,0){1}}
   \put( 10, 1){\line(0,1){24}}
   \put( 10, 1){\circle*{2}}
}

\newcommand{\cicard}[3]{
\setlength{\unitlength}{1.6pt}
\begin{picture}(20,70)(0,-16)
   \thinlines
   \multiput(  0, 0)(20,0){2}{\line(0,1){50}}
   \multiput(  0, 0)(0,50){2}{\line(1,0){20}}
   \put(  0, 50){\makebox(20, 7){{\small $#1$}}}
   \put(  0, -9){\makebox(20, 8){{\small $#2$}}}
   \put(  0,-16){\makebox(20, 6){{\small $#3$}}}
   \thicklines
  \stoa
  \atob
  \btoc
   \put(10,32){\makebox(0,0){$.$}}
   \put(10,29){\makebox(0,0){$.$}}
   \put(10,26){\makebox(0,0){$.$}}
  \dtoe
\end{picture}
}
\newcommand{\aijcard}[3]{
\setlength{\unitlength}{1.6pt}
\begin{picture}(20,70)(0,-16)
   \thinlines
   \multiput(  0, 0)(20,0){2}{\line(0,1){50}}
   \multiput(  0, 0)(0,50){2}{\line(1,0){20}}
   \put(  0, 50){\makebox(20, 7){{\small $#1$}}}
   \put(  0, -9){\makebox(20, 8){{\small $#2$}}}
   \put(  0,-16){\makebox(20, 6){{\small $#3$}}}
   \thicklines
   \atoa
   \put(15, 16){\makebox(0,0){$ \vdots $}}
   \btob
   \ctos
   \put( -6,26){\makebox(0,0){$\scriptstyle{j \text{th}}$}}
   \dtoc
   \put(15, 33){\makebox(0,0){$ \vdots $}}
   \etod
\end{picture}
}
\newcommand{\nijcard}[3]{
\setlength{\unitlength}{1.6pt}
\begin{picture}(20,70)(0,-16)
   \thinlines
   \multiput(  0, 0)(20,0){2}{\line(0,1){50}}
   \multiput(  0, 0)(0,50){2}{\line(1,0){20}}
   \put(  0, 50){\makebox(20, 7){{\small $#1$}}}
   \put(  0, -9){\makebox(20, 8){{\small $#2$}}}
   \put(  0,-16){\makebox(20, 6){{\small $#3$}}}
   \thicklines
   \stoa
   \atob
   \put(15,23){\makebox(0,0){$ \vdots $}}
   \btoc
   \ctos
   \put( -6,26){\makebox(0,0){$\scriptstyle{j \text{th}}$}}
   \dtod
   \put(15, 40){\makebox(0,0){$ \vdots $}}
   \etoe
\end{picture}
}
\newcommand{\kicard}[3]{
\setlength{\unitlength}{1.6pt}
\begin{picture}(20,70)(0,-16)
   \thinlines
   \multiput(  0, 0)(20,0){2}{\line(0,1){50}}
   \multiput(  0, 0)(0,50){2}{\line(1,0){20}}
   \put(  0, 50){\makebox(20, 7){{\small $#1$}}}
   \put(  0, -9){\makebox(20, 8){{\small $#2$}}}
   \put(  0,-16){\makebox(20, 6){{\small $#3$}}}
   \thicklines
   \stos
   \atoa
   \btob
   \put(10, 27){\makebox(0,0){$ \vdots $}}
   \dtod
\end{picture}
}
\newcommand{\aijcicard}{
\setlength{\unitlength}{1.6pt}
\begin{picture}(20,70)(0,-16)
   \thinlines
   \put(  0, 0){\line(0,1){50}}
   \multiput( 20,6)(0, 4){11}{\line(0,1){2}}
   \multiput(  0, 0)(0,50){2}{\line(1,0){20}}
   \put(  0, -9){\makebox(20, 8){$ A_i^{(j)}$}}
   \linethickness{1.5pt}
   \put( 10, 1){\line(1,0){10}}
   \put( 12, 3){$\scriptstyle{\text{identify }}$}
   \thicklines
   \atoa
   \put(15, 16){\makebox(0,0){$ \vdots $}}
   \btob
   \ctos
   \put( -7,26){\makebox(0,0){$\scriptstyle{j \text{th}}$}}
   \dtoc
   \put(15, 33){\makebox(0,0){$ \vdots $}}
   \etod
\end{picture}
\setlength{\unitlength}{1.6pt}
\begin{picture}(20,70)(0,-16)
   \thinlines
   \put( 20, 0){\line(0,1){50}}
   \multiput(  0, 0)(0,50){2}{\line(1,0){20}}
   \put(  0, -9){\makebox(20, 8){$ C_{i-1}$}}
   \linethickness{1.5pt}
   \put( 0, 1){\line(1,0){10}}
   \thicklines
   \stoa
   \atob
   \put(15,23){\makebox(0,0){$ \vdots $}}
   \btoc
   \ctod
   \put(15,39){\makebox(0,0){$ \vdots $}}
   \dtoe
\end{picture}
}
\newcommand{\crzer}[3]{
\setlength{\unitlength}{1.6pt}
\begin{picture}(20,70)(0,-16)
   \thinlines
   \multiput(  0, 0)(20,0){2}{\line(0,1){50}}
   \multiput(  0, 0)(0,50){2}{\line(1,0){20}}
   \put(  0, 50){\makebox(20, 7){{\small $#1$}}}
   \put(  0, -9){\makebox(20, 8){{\small $#2$}}}
   \put(  0,-16){\makebox(20, 6){{\small $#3$}}}
   \thicklines
   \stoa
\end{picture}
}
\newcommand{\crone}[3]{
\setlength{\unitlength}{1.6pt}
\begin{picture}(20,70)(0,-16)
   \thinlines
   \multiput(  0, 0)(20,0){2}{\line(0,1){50}}
   \multiput(  0, 0)(0,50){2}{\line(1,0){20}}
   \put(  0, 50){\makebox(20, 7){{\small $#1$}}}
   \put(  0, -9){\makebox(20, 8){{\small $#2$}}}
   \put(  0,-16){\makebox(20, 6){{\small $#3$}}}
   \thicklines
   \stoa
   \atob
\end{picture}
}
\newcommand{\crtwo}[3]{
\setlength{\unitlength}{1.6pt}
\begin{picture}(20,70)(0,-16)
   \thinlines
   \multiput(  0, 0)(20,0){2}{\line(0,1){50}}
   \multiput(  0, 0)(0,50){2}{\line(1,0){20}}
   \put(  0, 50){\makebox(20, 7){{\small $#1$}}}
   \put(  0, -9){\makebox(20, 8){{\small $#2$}}}
   \put(  0,-16){\makebox(20, 6){{\small $#3$}}}
   \thicklines
   \stoa
   \atob
   \btoc
\end{picture}
}
\newcommand{\crthr}[3]{
\setlength{\unitlength}{1.6pt}
\begin{picture}(20,70)(0,-16)
   \thinlines
   \multiput(  0, 0)(20,0){2}{\line(0,1){50}}
   \multiput(  0, 0)(0,50){2}{\line(1,0){20}}
   \put(  0, 50){\makebox(20, 7){{\small $#1$}}}
   \put(  0, -9){\makebox(20, 8){{\small $#2$}}}
   \put(  0,-16){\makebox(20, 6){{\small $#3$}}}
   \thicklines
   \stoa
   \atob
   \btoc
   \ctod
\end{picture}
}
\newcommand{\anoneone}[3]{
\setlength{\unitlength}{1.6pt}
\begin{picture}(20,70)(0,-16)
   \thinlines
   \multiput(  0, 0)(20,0){2}{\line(0,1){50}}
   \multiput(  0, 0)(0,50){2}{\line(1,0){20}}
   \put(  0, 50){\makebox(20, 7){{\small $#1$}}}
   \put(  0, -9){\makebox(20, 8){{\small $#2$}}}
   \put(  0,-16){\makebox(20, 6){{\small $#3$}}}
   \thicklines
   \atos
\end{picture}
}
\newcommand{\antwoone}[3]{
\setlength{\unitlength}{1.6pt}
\begin{picture}(20,70)(0,-16)
   \thinlines
   \multiput(  0, 0)(20,0){2}{\line(0,1){50}}
   \multiput(  0, 0)(0,50){2}{\line(1,0){20}}
   \put(  0, 50){\makebox(20, 7){{\small $#1$}}}
   \put(  0, -9){\makebox(20, 8){{\small $#2$}}}
   \put(  0,-16){\makebox(20, 6){{\small $#3$}}}
   \thicklines
   \atos
   \btoa
\end{picture}
}
\newcommand{\antwotwo}[3]{
\setlength{\unitlength}{1.6pt}
\begin{picture}(20,70)(0,-16)
   \thinlines
   \multiput(  0, 0)(20,0){2}{\line(0,1){50}}
   \multiput(  0, 0)(0,50){2}{\line(1,0){20}}
   \put(  0, 50){\makebox(20, 7){{\small $#1$}}}
   \put(  0, -9){\makebox(20, 8){{\small $#2$}}}
   \put(  0,-16){\makebox(20, 6){{\small $#3$}}}
   \thicklines
   \btos
   \atoa
\end{picture}
} 
\newcommand{\anthrone}[3]{
\setlength{\unitlength}{1.6pt}
\begin{picture}(20,70)(0,-16)
   \thinlines
   \multiput(  0, 0)(20,0){2}{\line(0,1){50}}
   \multiput(  0, 0)(0,50){2}{\line(1,0){20}}
   \put(  0, 50){\makebox(20, 7){{\small $#1$}}}
   \put(  0, -9){\makebox(20, 8){{\small $#2$}}}
   \put(  0,-16){\makebox(20, 6){{\small $#3$}}}
   \thicklines
   \atos
   \btoa
   \ctob
\end{picture}
}
\newcommand{\anthrtwo}[3]{
\setlength{\unitlength}{1.6pt}
\begin{picture}(20,70)(0,-16)
   \thinlines
   \multiput(  0, 0)(20,0){2}{\line(0,1){50}}
   \multiput(  0, 0)(0,50){2}{\line(1,0){20}}
   \put(  0, 50){\makebox(20, 7){{\small $#1$}}}
   \put(  0, -9){\makebox(20, 8){{\small $#2$}}}
   \put(  0,-16){\makebox(20, 6){{\small $#3$}}}
   \thicklines
   \ctob
   \btos
   \atoa
\end{picture}
}
\newcommand{\anthrthr}[3]{
\setlength{\unitlength}{1.6pt}
\begin{picture}(20,70)(0,-16)
   \thinlines
   \multiput(  0, 0)(20,0){2}{\line(0,1){50}}
   \multiput(  0, 0)(0,50){2}{\line(1,0){20}}
   \put(  0, 50){\makebox(20, 7){{\small $#1$}}}
   \put(  0, -9){\makebox(20, 8){{\small $#2$}}}
   \put(  0,-16){\makebox(20, 6){{\small $#3$}}}
   \thicklines
   \ctos
   \btob
   \atoa
\end{picture}
}
\newcommand{\anforone}[3]{
\setlength{\unitlength}{1.6pt}
\begin{picture}(20,70)(0,-16)
   \thinlines
   \multiput(  0, 0)(20,0){2}{\line(0,1){50}}
   \multiput(  0, 0)(0,50){2}{\line(1,0){20}}
   \put(  0, 50){\makebox(20, 7){{\small $#1$}}}
   \put(  0, -9){\makebox(20, 8){{\small $#2$}}}
   \put(  0,-16){\makebox(20, 6){{\small $#3$}}}
   \thicklines
   \atos
   \btoa
   \ctob
   \dtoc
\end{picture}
}
\newcommand{\anfortwo}[3]{
\setlength{\unitlength}{1.6pt}
\begin{picture}(20,70)(0,-16)
   \thinlines
   \multiput(  0, 0)(20,0){2}{\line(0,1){50}}
   \multiput(  0, 0)(0,50){2}{\line(1,0){20}}
   \put(  0, 50){\makebox(20, 7){{\small $#1$}}}
   \put(  0, -9){\makebox(20, 8){{\small $#2$}}}
   \put(  0,-16){\makebox(20, 6){{\small $#3$}}}
   \thicklines
   \atoa
   \btos
   \ctob
   \dtoc
\end{picture}
}
\newcommand{\anforthr}[3]{
\setlength{\unitlength}{1.6pt}
\begin{picture}(20,70)(0,-16)
   \thinlines
   \multiput(  0, 0)(20,0){2}{\line(0,1){50}}
   \multiput(  0, 0)(0,50){2}{\line(1,0){20}}
   \put(  0, 50){\makebox(20, 7){{\small $#1$}}}
   \put(  0, -9){\makebox(20, 8){{\small $#2$}}}
   \put(  0,-16){\makebox(20, 6){{\small $#3$}}}
   \thicklines
   \atoa
   \btob
   \ctos
   \dtoc
\end{picture}
}
\newcommand{\anforfor}[3]{
\setlength{\unitlength}{1.6pt}
\begin{picture}(20,70)(0,-16)
   \thinlines
   \multiput(  0, 0)(20,0){2}{\line(0,1){50}}
   \multiput(  0, 0)(0,50){2}{\line(1,0){20}}
   \put(  0, 50){\makebox(20, 7){{\small $#1$}}}
   \put(  0, -9){\makebox(20, 8){{\small $#2$}}}
   \put(  0,-16){\makebox(20, 6){{\small $#3$}}}
   \thicklines
   \atoa
   \btob
   \ctoc
   \dtos
\end{picture}
}
\newcommand{\nuoneone}[3]{
\setlength{\unitlength}{1.6pt}
\begin{picture}(20,70)(0,-16)
   \thinlines
   \multiput(  0, 0)(20,0){2}{\line(0,1){50}}
   \multiput(  0, 0)(0,50){2}{\line(1,0){20}}
   \put(  0, 50){\makebox(20, 7){{\small $#1$}}}
   \put(  0, -9){\makebox(20, 8){{\small $#2$}}}
   \put(  0,-16){\makebox(20, 6){{\small $#3$}}}
   \thicklines
   \atos
   \stoa
\end{picture}
}
\newcommand{\nutwoone}[3]{
\setlength{\unitlength}{1.6pt}
\begin{picture}(20,70)(0,-16)
   \thinlines
   \multiput(  0, 0)(20,0){2}{\line(0,1){50}}
   \multiput(  0, 0)(0,50){2}{\line(1,0){20}}
   \put(  0, 50){\makebox(20, 7){{\small $#1$}}}
   \put(  0, -9){\makebox(20, 8){{\small $#2$}}}
   \put(  0,-16){\makebox(20, 6){{\small $#3$}}}
   \thicklines
   \atos
   \stoa
   \btob
\end{picture}
}
\newcommand{\nutwotwo}[3]{
\setlength{\unitlength}{1.6pt}
\begin{picture}(20,70)(0,-16)
   \thinlines
   \multiput(  0, 0)(20,0){2}{\line(0,1){50}}
   \multiput(  0, 0)(0,50){2}{\line(1,0){20}}
   \put(  0, 50){\makebox(20, 7){{\small $#1$}}}
   \put(  0, -9){\makebox(20, 8){{\small $#2$}}}
   \put(  0,-16){\makebox(20, 6){{\small $#3$}}}
   \thicklines
   \stoa
   \atob
   \btos
\end{picture}
}
\newcommand{\nuthrone}[3]{
\setlength{\unitlength}{1.6pt}
\begin{picture}(20,70)(0,-16)
   \thinlines
   \multiput(  0, 0)(20,0){2}{\line(0,1){50}}
   \multiput(  0, 0)(0,50){2}{\line(1,0){20}}
   \put(  0, 50){\makebox(20, 7){{\small $#1$}}}
   \put(  0, -9){\makebox(20, 8){{\small $#2$}}}
   \put(  0,-16){\makebox(20, 6){{\small $#3$}}}
   \thicklines
   \atos
   \stoa
   \btob
   \ctoc
\end{picture}
}
\newcommand{\nuthrtwo}[3]{
\setlength{\unitlength}{1.6pt}
\begin{picture}(20,70)(0,-16)
   \thinlines
   \multiput(  0, 0)(20,0){2}{\line(0,1){50}}
   \multiput(  0, 0)(0,50){2}{\line(1,0){20}}
   \put(  0, 50){\makebox(20, 7){{\small $#1$}}}
   \put(  0, -9){\makebox(20, 8){{\small $#2$}}}
   \put(  0,-16){\makebox(20, 6){{\small $#3$}}}
   \thicklines
   \stoa
   \atob
   \btos
   \ctoc
\end{picture}
}
\newcommand{\nuthrthr}[3]{
\setlength{\unitlength}{1.6pt}
\begin{picture}(20,70)(0,-16)
   \thinlines
   \multiput(  0, 0)(20,0){2}{\line(0,1){50}}
   \multiput(  0, 0)(0,50){2}{\line(1,0){20}}
   \put(  0, 50){\makebox(20, 7){{\small $#1$}}}
   \put(  0, -9){\makebox(20, 8){{\small $#2$}}}
   \put(  0,-16){\makebox(20, 6){{\small $#3$}}}
   \thicklines
   \stoa
   \atob
   \btoc
   \ctos
\end{picture}
}
\newcommand{\nuforone}[3]{
\setlength{\unitlength}{1.6pt}
\begin{picture}(20,70)(0,-16)
   \thinlines
   \multiput(  0, 0)(20,0){2}{\line(0,1){50}}
   \multiput(  0, 0)(0,50){2}{\line(1,0){20}}
   \put(  0, 50){\makebox(20, 7){{\small $#1$}}}
   \put(  0, -9){\makebox(20, 8){{\small $#2$}}}
   \put(  0,-16){\makebox(20, 6){{\small $#3$}}}
   \thicklines
   \atos
   \stoa
   \btob
   \ctoc
   \dtod
\end{picture}
}
\newcommand{\nufortwo}[3]{
\setlength{\unitlength}{1.6pt}
\begin{picture}(20,70)(0,-16)
   \thinlines
   \multiput(  0, 0)(20,0){2}{\line(0,1){50}}
   \multiput(  0, 0)(0,50){2}{\line(1,0){20}}
   \put(  0, 50){\makebox(20, 7){{\small $#1$}}}
   \put(  0, -9){\makebox(20, 8){{\small $#2$}}}
   \put(  0,-16){\makebox(20, 6){{\small $#3$}}}
   \thicklines
   \stoa
   \atob
   \btos
   \ctoc
   \dtod
\end{picture}
}
\newcommand{\nuforthr}[3]{
\setlength{\unitlength}{1.6pt}
\begin{picture}(20,70)(0,-16)
   \thinlines
   \multiput(  0, 0)(20,0){2}{\line(0,1){50}}
   \multiput(  0, 0)(0,50){2}{\line(1,0){20}}
   \put(  0, 50){\makebox(20, 7){{\small $#1$}}}
   \put(  0, -9){\makebox(20, 8){{\small $#2$}}}
   \put(  0,-16){\makebox(20, 6){{\small $#3$}}}
   \thicklines
   \stoa
   \atob
   \btoc
   \ctos
   \dtod
\end{picture}
}
\newcommand{\nuforfor}[3]{
\setlength{\unitlength}{1.6pt}
\begin{picture}(20,70)(0,-16)
   \thinlines
   \multiput(  0, 0)(20,0){2}{\line(0,1){50}}
   \multiput(  0, 0)(0,50){2}{\line(1,0){20}}
   \put(  0, 50){\makebox(20, 7){{\small $#1$}}}
   \put(  0, -9){\makebox(20, 8){{\small $#2$}}}
   \put(  0,-16){\makebox(20, 6){{\small $#3$}}}
   \thicklines
   \stoa
   \atob
   \btoc
   \ctod
   \dtos
\end{picture}
}
\newcommand{\kozer}[3]{
\setlength{\unitlength}{1.6pt}
\begin{picture}(20,70)(0,-16)
   \thinlines
   \multiput(  0, 0)(20,0){2}{\line(0,1){50}}
   \multiput(  0, 0)(0,50){2}{\line(1,0){20}}
   \put(  0, 50){\makebox(20, 7){{\small $#1$}}}
   \put(  0, -9){\makebox(20, 8){{\small $#2$}}}
   \put(  0,-16){\makebox(20, 6){{\small $#3$}}}
   \thicklines
   \put( 10, 0){\line(0,1){6}}
   \put( 10, 1){\circle*{2}}
\end{picture}
}
\newcommand{\koone}[3]{
\setlength{\unitlength}{1.6pt}
\begin{picture}(20,70)(0,-16)
   \thinlines
   \multiput(  0, 0)(20,0){2}{\line(0,1){50}}
   \multiput(  0, 0)(0,50){2}{\line(1,0){20}}
   \put(  0, 50){\makebox(20, 7){{\small $#1$}}}
   \put(  0, -9){\makebox(20, 8){{\small $#2$}}}
   \put(  0,-16){\makebox(20, 6){{\small $#3$}}}
   \thicklines
   \stos
   \atoa
\end{picture}
}
\newcommand{\kotwo}[3]{
\setlength{\unitlength}{1.6pt}
\begin{picture}(20,70)(0,-16)
   \thinlines
   \multiput(  0, 0)(20,0){2}{\line(0,1){50}}
   \multiput(  0, 0)(0,50){2}{\line(1,0){20}}
   \put(  0, 50){\makebox(20, 7){{\small $#1$}}}
   \put(  0, -9){\makebox(20, 8){{\small $#2$}}}
   \put(  0,-16){\makebox(20, 6){{\small $#3$}}}
   \thicklines
   \stos
   \atoa
   \btob
\end{picture}
}
\newcommand{\kothr}[3]{
\setlength{\unitlength}{1.6pt}
\begin{picture}(20,70)(0,-16)
   \thinlines
   \multiput(  0, 0)(20,0){2}{\line(0,1){50}}
   \multiput(  0, 0)(0,50){2}{\line(1,0){20}}
   \put(  0, 50){\makebox(20, 7){{\small $#1$}}}
   \put(  0, -9){\makebox(20, 8){{\small $#2$}}}
   \put(  0,-16){\makebox(20, 6){{\small $#3$}}}
   \thicklines
   \stos
   \atoa
   \btob
   \ctoc
\end{picture}
}
\newcommand{\dumcdots}{
\setlength{\unitlength}{1.6pt}
\begin{picture}(10,58)(0,-16)
   \put(0, 18){\makebox(10, 6){\Large{$\cdots$}}}
\end{picture}
}
\newcommand{\dumcomma}{
\setlength{\unitlength}{1.6pt}
\begin{picture}(10,58)(0,-16)
   \put(0, 0){\makebox(5,5){{\Large{,}}}}
\end{picture}
}

\newcommand{\dumarrow}{
\setlength{\unitlength}{1.6pt}
\begin{picture}(10,58)(0,-16)
   \put(0, 18){\makebox(10, 6){\Large{$\Longrightarrow$}}}
\end{picture}
}
\newcommand{\rsideilong}[1]{
\setlength{\unitlength}{1.6pt}
\begin{picture}(10,70)(0,-16)
   \put(0,28){{$\left. \makebox(0,18){} \right\} \scriptstyle{#1} $ }}
\end{picture}
}
\newcommand{\rsideihi}[1]{
\setlength{\unitlength}{1.6pt}
\begin{picture}(10,70)(0,-16)
   \put(0,28){{$\left. \makebox(0,6){} \right\} \scriptstyle{#1} $ }}
\end{picture}
}
\newcommand{\rsideiji}[2]{
\setlength{\unitlength}{1.6pt}
\begin{picture}(10,70)(0,-16)
   \put(0,36){{$\left. \makebox(0,6){} \right\} \scriptstyle{#1} $ }}
   \put(0,20){{$\left. \makebox(0,6){} \right\} \scriptstyle{#2} $ }}
\end{picture}
}
\newcommand{\rsideilobig}[1]{
\setlength{\unitlength}{1.6pt}
\begin{picture}(10,70)(0,-16)
   \put(0,20){{$\left. \makebox(0,15){} \right\} \scriptstyle{#1} $ }}
\end{picture}
}
\newcommand{\upze}{
\setlength{\unitlength}{1.6pt}
\begin{picture}(20,70)(0,-16)
   \thinlines
   \multiput(  0, 0)(0,50){1}{\line(1,0){20}}
   \put(  0, -9){\makebox(20, 8){$a$}}
   \qbezier(10,32)(10,42)(20,42)
   \put( 10, 1){\line(0,1){32}}
   \put( 10, 1){\circle*{2}}

   \put(  1, 1){\oval(18,18)[tr]}
\end{picture}
}

\newcommand{\upzc}{
\setlength{\unitlength}{1.6pt}
\begin{picture}(20,70)(0,-16)
   \thinlines
   \multiput(  0, 0)(0,50){1}{\line(1,0){20}}
   \put(  0, -9){\makebox(20, 8){$b$}}

   \qbezier(10,16)(10,26)(20,26)
   \put( 10, 1){\line(0,1){16}}
   \put( 10, 1){\circle*{2}}
   \put( 0, 42){\line(1,0){20}}
\end{picture}
}

\newcommand{\nestcard}{
\setlength{\unitlength}{1.6pt}
\begin{picture}(20,70)(0,-16)
   \thinlines
   \multiput(  0, 0)(20,0){2}{\line(0,1){50}}
   \multiput(  0, 0)(0,50){2}{\line(1,0){20}}
   \put(  0, -9){\makebox(20, 8){$c$}}
   \linethickness{1pt}
   \qbezier( 0,26)(10,26)(10,16)
   \put(  0,26){\line(1,0){1}}
   \put( 10, 1){\line(0,1){16}}
   \put( 10, 1){\circle*{2}}
   \etod{}
\end{picture}
}
\newcommand{\dndz}{
\setlength{\unitlength}{1.6pt}
\begin{picture}(20,70)(0,-16)
   \thinlines
   \multiput(  0, 0)(0,50){1}{\line(1,0){20}}
   \put(  0, -9){\makebox(20, 8){$d$}}
   \qbezier( 0,34)(10,34)(10,24)
   \put(  0,34){\line(1,0){1}}
   \put( 10, 1){\line(0,1){24}}
   \put( 10, 1){\circle*{2}}
   \put( 19, 1){\oval(18,18)[tl]}
\end{picture}
}

\newcommand{\crosscard}{
\setlength{\unitlength}{1.6pt}
\begin{picture}(20,70)(0,-16)
  \thinlines
   \multiput(  0, 0)(20,0){2}{\line(0,1){50}}
   \multiput(  0, 0)(0,50){2}{\line(1,0){20}}
   \put(  0, -9){\makebox(20, 8){$c$}}
   \linethickness{1pt}
   \qbezier( 0,42)(10,42)(10,32)
   \put(  0,42){\line(1,0){1}}
   \put( 10, 1){\line(0,1){32}}
   \put( 10, 1){\circle*{2}}
   \put( 0, 26){\line(1,0){20}}
\end{picture}
}
\newcommand{\dncz}{
\setlength{\unitlength}{1.6pt}
\begin{picture}(20,70)(0,-16)
   \thinlines
   \multiput(  0, 0)(0,50){1}{\line(1,0){20}}
   \put(  0, -9){\makebox(20, 8){$d$}}
   \qbezier( 0,26)(10,26)(10,16)
   \put(  0,26){\line(1,0){1}}
   \put( 10, 1){\line(0,1){16}}
   \put( 10, 1){\circle*{2}}
   \put( 19, 1){\oval(18,18)[tl]}
\end{picture}
}

\begin{document}


\begin{center}
{\bf \Large A Poisson Type Operator Deformed by Generalized Fibonacci 
Numbers and Its Combinatorial Moment Formula\footnote{Submitted for publication on July 29, 2025.}}

\bigskip

{
         Nobuhiro ASAI\footnote{
         Department of Mathematics,
         Aichi University of Education, 
       Kariya 448-8542, Japan.
 		{\tt nasai@auecc.aichi-edu.ac.jp}}, 
\vspace{0.5mm}
		Marek BO\.ZEJKO\footnote{
		Institute of Mathematics, 
		University of Wroc{\l}aw, 
		Plac. Grunwaldzki 2, 50-384, Wroc{\l}aw, Poland.
		{\tt bozejko@gmail.com}}, \\
\vspace{0.5mm}
    	Lahcen OUSSI\footnote{
		Department of Mathematics and Cybernetics, 
		Faculty of Economics and Finance, 
		Wroc{\l}aw University of Economics and Business,
		 ul. Komandorska 118/120, 53-345, Wroc{\l}aw, Poland.
		{\tt lahcen.oussi@ue.wroc.pl}, {\tt oussimaths@gmail.com}}, 
		and 
         Hiroaki YOSHIDA\footnote{
         Department of Mathematics and Information Science, 
		Faculty of Science, Josai University, 
	Tokyo 102‑0093, Japan, \\
    {\tt hyoshida@josai.ac.jp},  
	and Ochanomizu University, 
	Tokyo 112-8610, Japan,
	{\tt yoshida@is.ocha.ac.jp}}
 }

\medskip
July 25,  2025
\end{center}

\begin{abstract}
We introduce a two-parameter deformation of the classical Poisson distribution 
from the viewpoint of noncommutative probability theory, 
by defining a $(q,t)$-Poisson type operator (random variable) 
on the $(q,t)$-Fock space \cite{Bl12}  (See also \cite{BY06, AY20}).  
From the analogous viewpoint of the classical Poisson limit theorem
in probability theory, 
we are naturally led to a family of orthogonal polynomials, 
which we call the $(q,t)$-Charlier polynomials.
These generalize the $q$-Charlier polynomials of 
Saitoh-Yoshida \cite{SY00a, SY00b}
and reflect deeper combinatorial symmetries 
through the additional deformation parameter $t$.
A central feature of this paper is the derivation of 
a combinatorial moment formula of the $(q,t)$-Poisson type operator
and the $(q,t)$-Poisson distribution.
This is accomplished 
by means of a card arrangement technique, 
which encodes set partitions together with crossing and nesting statistics. 
The resulting expression naturally exhibits a duality between these statistics,
arising from a structure rooted in generalized Fibonacci numbers.
Our approach provides a concrete framework 
where methods in combinatorics and theory of orthogonal polynomials 
are used to investigate the probabilistic properties arising 
from the $(q,t)$-deformation.

\medskip
\noindent
{\bf Keywords:} 
deformation,
Poisson type operator (random variable), 
Charlier polynomials, Poisson distribution,
set partitions, partition statistics, restricted crossing and nesting,
generalized Fibonacci number, 
combinatorial duality,  
combinatorial moment formula

\smallskip 

\noindent
{\bf 2020 Mathematics Subject Classification:} 
46L53, 05A17, 33D45, 60E99, 05A30

\end{abstract}
\section{Introduction}
A {\it noncommutative} (or {\it quantum}) probability space is a unital 
(possibly noncommutative) algebra $\mathcal{A}$ together with a linear 
functional 
$\varphi : \mathcal{A} \to \mathbb{C}$, such that $\varphi(1) = 1$.
If $\mathcal{A}$ is a $C^*$-algebra and $\varphi$ is a state, then 
$\big( \mathcal{A}, \varphi \big)$ is 
called a {\it $C^*$-probability space}.
An operator in $\mathcal{A}$ is regarded as a {\it noncommutative 
random variable} and the 
{\it distribution of $x \in \mathcal{A}$} with respect to $\varphi$ is 
determined by the 
linear functional $\mu$ on $\mathbb{C}[X]$ (the polynomials in one variable) 
by 
$$
    \mu : \mathbb{C}[X] \ni P \longmapsto \varphi \big( P(X) \big) \in \mathbb{C}.
$$
Considered in the $C^*$-probability context, 
the distribution $\mu$ of a 
self-adjoint operator $x \in \mathcal{A}$ can be extended to, and 
identified with the (compactly supported) probability measure 
$\mu$ on $\mathbb{R}$ by 
$$
  \varphi \big( P(X) \big) = \int_{\mathbb{R}} P(t) \, \text{d} \mu(t), 
  \quad P \in {\mathbb{C}}[X].
$$

The following question has been considered in the literature, what
probability distribution will be obtained 
in a noncommutative central limit.
That is, what is the Gaussian counterpart in noncommutative setting ?
In noncommutative probabilistic framework regarding the  Boson Fock space 
by Hudson-Parthasarathy \cite{HP84}, 
it is the fundamental fact that 
the Gaussian operator (field operator), that is,  a sum of the creation $A^{\dagger}$ 
and the annihilation $A$ operators on Boson Fock space gives rise to 
Gaussian distribution in the vacuum state.
It is well-known that 
$A^{\dagger}$ and $A$ satisfy the canonical commutation relation (CCR).
In \cite{BS91}, Bo{\.z}ejko-Speicher introduced the $q$-Fock space
to the field of noncommutative probability theory
by deforming the inner product of free (full) Fock space 
with the positive definite function on the symmetric group. 
This gives an interpolation between the Boson (symmetric) and 
the Fermion (anti-symmetric) Fock spaces and, especially, the case 
$q=0$ of which yields the canonical model in the free probability 
theory (See \cite{VDN92}, for instance).
The corresponding $q$-creation and annihilation operators $A_q^{\dagger}$ 
and $A_q$, respectively, satisfy the commutation relation with 
the $q$-commutator ($q$-CCR), 
$A_q \, A_q^{\dagger} - q \, A_q^{\dagger} \, A_q = {\bm 1}$,
and the $q$-Gaussian operator ($q$-field operator) 
$A_q^{\dagger} + A_q$ gives rise to the $q$-Gaussian 
distribution, which is the orthogonalizing probability measures for 
the Rogers $q$-Hermite polynomials under the appropriate 
rescaling.   

On the other hand, the Poisson distribution also plays a fundamental role in 
probability theory.  In particular, using Gaussian and Poisson distributions, 
one can construct arbitrary infinitely divisible distributions via 
L{\'e}vy-Khintchine representation. 
Its noncommutative probabilistic setting, the $q$-Poisson operator (random variable) 
has been introduced on the $q$-Fock space 
and the $q$-Poisson distribution is known to be orthogonalized by 
the $q$-Charlier polynomials.  See papers \cite{SY00a, SY00b} by Saitoh-Yoshida.

In this paper, we investigate a two parameter deformation of 
the $q$-Poisson model  
by incorporating an additional deformation parameter $t$.
This leads to a new family of orthogonal polynomials, 
which we call the $(q,t)$-Charlier polynomials.
These are naturally defined from the viewpoint of 
the classical Poisson limit theorem in probability theory
and reflect the algebraic structure of underlying 
the $(q,t)$-Fock space ${\mathcal F}_{q,t}$ constructed by Blitvić \cite{Bl12}
(See also  \cite{BY06, AY20}).
We remark that for $q = 0$ and $t \in (0,1]$, 
this framework naturally includes, to the best of our knowledge, 
the orthogonal polynomials associated with 
what we call the \emph{$t$-free Poisson (Marchenko--Pastur)} distribution, 
which serves as the counterpart of the $t$-free Gaussian (semicircle) distribution discussed in \cite{Bl12}, 
and it extends the classical and free Poisson cases, 
which correspond to $q \to 1, t = 1$ and $q = 0, t = 1$, respectively.

Our $(q,t)$-Poisson type operator is realized as a linear combination of 
the the $(q,t)$-creation operator $A_{q,t}^{\dagger}$,
the $(q,t)$-annihilation operator $A_{q,t}$, $(q,t)$-number operator $N_{q,t}$ , and the scalar operator
acting on ${\mathcal F}_{q,t}$.   It is known \cite{Bl12} that $A_{q,t}^{\dagger}$ and $A_{q,t}$ 
satisfy the deformed commutation relation  
$A_{q,t} A_{q,t}^{\dagger} - q A_{q,t}^{\dagger} A_{q,t} = t^N$,
where $N$ is the number operator corresponding to the grading on ${\mathcal F}_{q,t}$.  
This commutation relation defines
the deformed quantum oscillator algebra of Chakrabarti and Jagannathan \cite{CJ91}.  
See also \cite{GK16, FKY22} and references therein.

A key result in this paper is a combinatorial moment formula 
for the $(q,t)$-Poisson type operator and the $(q,t)$-Poisson distribution. 
These moments are described in terms of set partitions 
and refined by partition statistics such as restricted crossings and nestings, 
which are inspired by the card arrangement technique for calculations of vacuum expectation. 
In particular, we derive moment expressions using a refined version 
of the {\it card arrangement technique}, which allows the encoding of operator products 
in a diagrammatic combinatorial form.
Our approach demonstrates how noncommutative 
probabilistic models can be analyzed through the viewpoint 
of combinatorics and orthogonal polynomials. 
In particular, it provides a concrete setting in which $(q,t)$-deformed combinatorial statistics, 
generalized Fibonacci numbers, and operator-theoretic constructions interact
to highlight and illustrate new probabilistic structures.

The paper is organized as follows. 
In Section \ref{sec2}, we introduce the 
$(q,t)$-Fock space, creation, and annihilation operators on this Fock space.
In Section \ref{sec:3}, 
we firstly explain 
an analogue of classical Poisson limit theorem from the 
point of the Jacobi parameters associated with 
Krawtchouk polynomials, quickly.
The recurrence formula for the orthogonal polynomials of 
the $(q,t)$-Poisson distribution is determined. 
Secondly,  
the $(q,t)$-Poisson type operator (random variable) will be 
defined. 
\if10
As we have mentioned above, in order to 
realize the deformed Poisson type operator (random variable), 
an appropriate gauge operator is required, 
and in our model the $(q,t)$-number operator 
is adopted at the gauge part. 
Hence,  our deformation in this paper can be  
regarded as a generalization of Saito-Yoshida type,
which is radically different from 
those of \cite{Ej20, AY25}.
\fi
\if10
A combinatorial representation of the moments of 
the $(q,t)$-Poisson distribution is given in this paper. 
\fi
In Section \ref{sec:set-partition}, we prepare combinatorial notions to state  
our moment formula in terms of the set partitions and their statistics.
In Section \ref{sec:moment}, we will present the combinatorial $n$-th moment formula 
of the $(q,t)$-Poisson distribution.  Our approach relies on the graphical 
interpretation of the  vacuum expectation of the product of operators based upon 
the card arrangement viewpoint.

%
\section{Preliminaries}\label{sec2}
In noncommutative probability theory, the Gaussian and Poisson processes 
can be realized on an appropriate Fock space 
by constructing related  creation and annihilation operators on that space. 
In the past years, probabilistic studies on Fock spaces have been 
done, for instance, in pioneer papers \cite{HP84, BS91, BKS97, M97, S97}, 
which provided the foundational spaces for the noncommutative probability theory. 
Related works, for example, \cite{BY06, Bl12, BEH15, AY20, O23, OW24, AY24, AY25},
aim at their generalization, refinement, and unification. 

This section is devoted to the construction of the $(q,t)$-Fock space 
of Blitvić \cite{Bl12} using the $\tau$-weighted sequence approach developed 
in \cite{BY06, AY20}.
In the next Section \ref{sec:3}, we discuss how this construction 
relates to our $(q,t)$-deformation of the Poisson distribution and 
Poisson type operator.

\subsection{Operators on Deformed Fock Space}
Let $\hilsp$ be a real Hilbert space equipped with the inner 
product $\langle \, \cdot \,  | \, \cdot \, \rangle$, and let
$\Omega$ be a distinguished unit vector, the so-called vacuum vector.
We denote by ${\mathcal F}_{\text{fin}} (\hilsp)$ the set of all the finite 
linear combinations of the elementary vectors 
$\xi_1 \otimes$ $\cdots$ $\otimes \xi_n$ $\in \hilsp^{\otimes n}$ \; 
$(n = 0, 1, 2, \ldots )$, where 
$\hilsp^{\otimes 0} = {\mathbb C} \Omega$ as convention.

The following notation of $(q,t)$-calculus should be prepared 
for later use.  
For $n \in {\mathbb N} \cup  \{ 0 \}$, let $[n]_{q,t}$ 
be the   {\it generalized Fibonacci number} given by
$$
  [n]_{q,t} := \frac{t^n - q^n}{t - q} = \sum_{k=1}^n t^{n-k}q^{k-1}, 
\ t\ne q, \ n\geq 1,
$$
and $[n]_{q,t}\to q^{n-1}n$ as $t\to q$ where $[0]_{q,t} = 0$ as convention.
$[n]_{q,1}$ is nothing but the $q$-number denoted by $[n]_q$ and hence $[n]_{1}=n$.
In this paper,  we refer to the {\it generalized Fibonacci number} as 
the {\it $(q,t)$-number} just for simplicity.
 For more details about the $(q,t)$-number and related topics, we refer the 
readers to \cite{GK16, FKY22, O22} 
and references therein.

Let us now recall the minimum about 
the $(q,t)$-Fock space
considered as the weighted $q/t$-Fock space 
with the weight $\tau_n = t^{n-1}$ $(n \ge 1)$
because of $[n]_{q,t}=t^{n-1}[n]_{q/t}$.
See \cite{BY06, AY20}. 
Since one can see easily $[n]_{q,t}=[n]_{t,q}$, 
let us assume $q\in (-1,1)$ and $|q|<t \leq 1$ without loss of generality,
and introduce the $(q,t)$-inner product    
$\big( \, \cdot  \, | \, \cdot \,  \big)_{q,t}$ 
on ${\mathcal F}_{\text{fin}} (\hilsp)$ by 
$$
  \big(  \xi_1 \otimes \cdots \otimes \xi_n \,  | \,  
          \eta_1 \otimes \cdots \otimes \eta_m  \big)_{q,t}
  =  \delta_{m, n} \, t^{\frac{n(n - 1)}{2}} \! \! 
      \sum_{\sigma \in {\mathfrak S}_n}  (q/t)^{i(\sigma)}
       \langle \xi_1 \, |  \, \eta_{\sigma(1)} \rangle \cdots 
       \langle \xi_n \, |  \, \eta_{\sigma(n)} \rangle ,
$$
where ${\mathfrak S}_n$ is the $n$-th symmetric group of 
permutations and $i(\sigma)$ is the number of inversions of the 
permutation $\sigma \in {\mathfrak S}_n$ defined by 
$$
 i(\sigma) = {\#} \big\{ \, (i,j) \, \big| \, 1 \le i < j \le n, \,  
                     \sigma(i) > \sigma (j) \big\}.
$$
The positivity of the inner product
 $\big(\, \cdot \,  | \, \cdot \, \big)_{q,t}$ 
is ensured in \cite{BS94, BY06, AY20}, which allows us to give the following 
definition:  

\begin{definition}
{\it The $(q,t)$-Fock space} ${\mathcal F}_{q,t} (\hilsp)$ is defined by
$$
	{\mathcal F}_{q,t} (\hilsp)
	= {\mathbb C} \Omega \oplus \bigoplus_{n=1}^{\infty}\hilsp^{\otimes n},
$$ 
the completion of ${\mathcal F}_{\text{fin}} (\hilsp)$ with respect 
to the inner product $\big( \, \cdot \,  | \, \cdot \, \big)_{q,t}$. 
\end{definition}

\begin{definition}
For a given $\xi \in \hilsp$, 
{\it the $(q,t)$-creation operator $A_{q,t}^{\dagger} (\xi)$} 
is defined by the canonical left creation,
\begin{align}\label{eq:a^*}
   &  A_{q,t}^{\dagger} (\xi) \, \Omega = \xi,  \notag \\
   & A_{q,t}^{\dagger} (\xi) \, (\xi_1 \otimes \cdots \otimes \xi_n) 
                    = \xi \otimes \xi_1 \otimes \cdots \otimes \xi_n, 
                  \quad n \geq 1.
\end{align}
{\it The $(q,t)$-annihilation operator $A_{q,t}(\xi)$}
 is defined by the adjoint operator of $A_{q,t}^{\dagger}(\xi)$ 
with respect to the inner product 
$\big(\, \cdot \,  | \, \cdot \, \big)_{q,t}$, 
that is, $A_{q,t}(\xi) = \Big( A_{q,t}^{\dagger} (\xi) \Big)^{*}$.
\end{definition}

By a direct consequence of the above definition,
the action of the $(q,t)$-annihilation operator on the elementary vectors 
has the form \cite{AY20},
\begin{equation}\label{eq:a-1}
	A_{q,t}(\xi) \, (\xi_1 \otimes \cdots \otimes \xi_n) 
	= t^{n-1}
	\displaystyle{
	\sum_{k=1}^n  
	\left(\frac{q}{t}\right)^{k-1} \langle \xi | \xi_k \rangle \, 
                     \xi_1 \otimes \cdots \otimes \stackrel{\vee}{\xi_k} 
                                   \otimes \cdots \otimes \xi_n },
      \quad n \geq 2,
\end{equation}
and hence one can get the following equivalent expression in \cite{Bl12}:

\begin{proposition}
 The $(q,t)$-annihilation operator $A_{q,t}(\xi)$ acts on the 
  elementary vectors as follows: 
\begin{align}\label{eq:a}
  &  A_{q,t}(\xi) \, \Omega = 0, \quad 
    A_{q,t}(\xi) \, \xi_1 = \langle \xi | \xi_{1} \rangle \, \Omega , \notag \\ 
  & A_{q,t}(\xi) \, (\xi_1 \otimes \cdots \otimes \xi_n) 
      = \displaystyle{\sum_{k=1}^n  
                     q^{k-1}t^{n-k} \langle \xi | \xi_k \rangle \, 
                     \xi_1 \otimes \cdots \otimes \stackrel{\vee}{\xi_k} 
                                   \otimes \cdots \otimes \xi_n },
      \quad n \ge 2,
\end{align}
 where $\stackrel{\vee}{\xi_k}$ means that $\xi_k$ should be deleted from 
 the tensor product.
\end{proposition}

Let us consider the state $\varphi$ for bounded operators on 
the $(q,t)$-Fock space ${\mathcal F}_{q,t} ({\mathscr H})$
given by the vacuum vector $\Omega$ as 
\begin{equation*}
   \varphi (b) = \left( b \, \Omega \, | \, \Omega \right)_{q,t},
   \qquad b \in \mathcal{B} \big( {\mathcal F}_{q,t} ({\mathscr H}) \big),
\end{equation*}
which is called {\it the vacuum expectation of $b$}.
One can employ 
$\big( \mathcal{B} \big( {\mathcal F}_{q,t} ({\mathscr H}) \big), 
     \varphi \big)$ 
as the noncommutative probability space, on which the model of the 
$(q,t)$-Poisson type operator (random variable) will be discussed.

It is known that the $(q,t)$-creation and the $(q,t)$-annihilation operators satisfy 
the following commutation relation \cite{Bl12}, 
\begin{equation}\label{eq:comrel}
    A_{q, t}(\xi) \, A_{q, t}^{\dagger} (\eta) 
       - q A_{q,t}^{\dagger}(\eta) \, A_{q,t}(\xi) 
     = \langle \xi | \eta \rangle \, t^N,
       \qquad \xi, \eta \in \hilsp,
\end{equation}
where the operator $t^N$ is defined on ${\mathcal F}_{q,t} ({\mathscr H})$ by 
\begin{equation*}
	t^N\Omega=\Omega, 
	\quad
	t^N(\xi_1\otimes \cdots \otimes \xi_n)
	=t^n\xi_1\otimes \cdots \otimes \xi_n, 
	\quad n\geq 1.
\end{equation*}

\begin{remark}
By the $\tau$-weight sequence approach \cite{BY06, AY20}, 
the commutation relation \eqref{eq:comrel} can be obtained
as a special choice of weight $\tau$.
See Section 3 of \cite{AY20} in detail.
\end{remark}

\begin{remark}
We note that the $(q,t)$- notation used in this paper 
differs from the $(q,s)$- and $(\alpha,q)$-deformation 
terminologies in \cite{AY25} and \cite{BEH15}, respectively.
\end{remark}

%
\section{Deformed Charlier Polynomials and Poisson Type Operator}\label{sec:3}
\subsection{Poisson as a limiting case of 
Binomial distribution}\label{sec:qt-Bin}

Before proceeding to the $(q,t)$-deformation, we recall a basic fact on orthogonal polynomials.
See \cite{Ch78, HO07} for details.

Let $\mu$ be a probability measure on ${\mathbb R}$ with finite moments of all orders. 
Then, it is well-known that there exists sequences of real numbers $\alpha_{n}\in {\mathbb R}$ 
and $\omega_{n}\geq 0$, so-called Jacobi parameters, 
such that the sequence of orthogonal polynomials $\{ P_{n}(x)\}$
with respect to $\mu$ is given by the recurrence relation,
\begin{align*}
	& P_{0}(x)=1, \ P_{1}(x)=x-\alpha_0, \\
	& P_{n+1}(x)=(x-\alpha_{n})P_{n}(x)-\omega_{n}P_{n-1}(x), \quad n\geq 1.
\end{align*}
In fact, the following orthogonality relation holds:
\begin{equation*}
	\int_{\mathbb R}P_{n}(x)P_{m}(x)d\mu(x)=\left(\prod_{i=1}^{n}\omega_{i}\right)\delta_{n, m}.
\end{equation*}
Conversely, the Farvard\rq{}s theorem \cite{Ch78} ensures the existence 
of a probability measure such that  for given 
parameters $\alpha_n$ and $\omega_n$, 
the sequence of polynomials determined by the 
above recurrence relation is orthogonal. 
Moreover, the probability measure is supported only in finitely many points 
if and only if there exists a number $m_0\geq 1$ such that $\omega_n=0$ for all $n\geq m_0$,
which means the sequence of polynomials is finite.

Let $\mu_{p}^{m}$ be the probability measure for the binomial  distribution $B(m, p)$,
\begin{equation}
	\mu_{p}^{m}(ds)=\sum_{k=0}^{m}\binom{m}{k}p^{k}(1-p)^{m-k}\delta_{k}(s),
	\quad 0<p<1,
\end{equation}
where $ds$ is the Lebesgue measure and $\delta_k(s)$ denotes the Dirac mass at $s=k$.
Orthogonal polynomials associated with $B(m,p)$ are so-called Krawtchouk polynomials 
determined by the Jacobi parameters \cite{Ch78}, 
\begin{equation}
	\alpha_{n}^{(m, p)}=mp+(1-2p)n, \quad\quad \omega_{n}^{(m, p)}=n(m-n+1)p(1-p).
\end{equation}
One can start a $(q,t)$-deformation of 
$\alpha_{n}^{(m, p)}$ and $\beta_{n}^{(m, p)}$ for $B(m,p)$.
By considering an $(q,t)$-analogue of the well-known Poisson limit theorem, 
one can consider a $(q,t)$-deformation of 
the Poisson distribution based on the idea of orthogonal polynomials.

\begin{definition}
For $q\in (-1,1)$ and $|q|<t\leq1$, 
the probability measure $(\mu_{p}^{m})_{q,t}$ induced from the Jacobi parameters, 
\begin{align}
	& (\alpha_{n}^{(m, p)})_{q,t}=mp+(1-2p)[n]_{q,t}, \label{alphalacobipq} \\
	& (\omega_{n}^{(m, p)})_{q,t}=[n]_{q,t}(m-[n-1]_{q,t})p(1-p)>0, \label{betalacobipq}
\end{align}
is called the $(q,t)$-deformed binomial distribution, denoted by 
$B_{q,t}(m, r)$, where $(\omega_{n}^{(m, r)})_{q,t}=0$ for all $n$ satisfying 
an inequality $[n-1]_{q,t}\geq m$. 
\end{definition}

Now, by taking the Poisson limits in the Jacobi parameters \eqref{alphalacobipq} 
and \eqref{betalacobipq}, that is, $m\to \infty$ 
and $p\to 0$ keeping the condition $mp=\lambda>0$, one can obtain 
\begin{align*}
	& (\alpha_{n}^{(m, p)})_{q,t}=mp+(1-2p)[n]_{q,t}
	\longrightarrow \lambda+[n]_{q,t},\\
	& (\omega_{n}^{(m, p)})_{q,t}=[n]_{q,t}(m-[n-1]_{q,t})p(1-p)
	 \longrightarrow \lambda[n]_{q,t}.		
\end{align*}

Therefore, as a direct consequence of considering the 
$(q,t)$-analogue of the classical Poisson limit theorem, 
we naturally arrive at the $(q,t)$-deformation of the Charlier polynomials
as follow.

\begin{definition}\label{def:qt-charlier}
(1) For $q\in (-1,1)$ and $|q|<t\leq1$ and $\lambda>0$, 
let $\mu_{q,t}^{\lambda}$ be  the orthogonalizing probability 
measure on $\mathbb R$ for the sequence of orthogonal polynomials 
$\{ C_n^{(q,t)} (\lambda ; x) \}$ determined by the following recurrence relation: 
 \begin{align}\label{eq:qt-op-c}
    & C_0^{(q,t)} (\lambda ; x) = 1, \quad C_1^{(q,t)} (\lambda ; x) = x - \lambda,  \notag\\
    & C_{n+1}^{(q,t)} (\lambda ; x) = 
        \left( x - (\lambda  + [n]_{q,t} ) \right) C_n^{(q,t)} (\lambda ; x) 
          - \lambda  [n]_{q,t} \, C_{n-1}^{(q,t)} (\lambda ; x), \quad n \ge 1. 
\end{align}

\noindent
(2) In this paper, $\mu_{q,t}^{\lambda}$ is called  the $(q,t)$-Poisson distribution
and  $\{C^{(q,t)}_n(\lambda;x)\}$ abbreviated as $\{C_{n}(x)\}$ from now on are 
called the $(q,t)$-Charlier polynomials.

\end{definition}

By the recurrence formula \eqref{eq:qt-op-c}, 
one can see that 
the first few terms of $(q,t)$-Charlier polynomials are 
given as follows:
\begin{align*}
	C_{1}(x) &=x-\lambda,\\
	C_{2}(x) &=x^2-(2\lambda+1)x+\lambda^2,\\
	C_{3}(x) &=x^3-(3\lambda+t+q+1)x^2+(3\lambda^2+(t+q)(\lambda+1)+\lambda)x-\lambda^3.
\end{align*} 

\begin{remark}
We remark that the $q$-deformation for $q\in [0,1)$ by Saitoh-Yoshida \cite{SY00a}
can be obtained as a special case if $t=1$.  
Since our deformation can be done for $q\in (-1,1)$ under 
$|q|<t\leq 1$, the deformation treated in this paper is considered as an 
important extension of the result in \cite{SY00a}. 
\end{remark}

\begin{remark}\label{rem:ej-poisson}
As one can see, $\{C_n(x)\}$ given in \eqref{eq:qt-op-c}
are different from $(q,t)$-Poisson polynomials $\{P_n^{(q,t)}(x)\}$ in \cite{Ej20} 
defined by 
\begin{align}\label{eq:ej-poisson}
	& P_0^{(q,t)} (x) = 1, \quad P_1^{(q,t)} (x) = x,  \notag\\
	& (x-[n]_{q,t})P_n^{(q,t)}(x)=P_{n+1}^{(q,t)}(x)+[n]_{q,t}P_{n-1}^{(q,t)}(x), \ n\geq 1.
\end{align}
\end{remark}

\subsection{$\bm{(q,t)}$-Poisson Type Operator}\label{sec:qt-poisson}
From now on, let us treat the $(q,t)$-Fock space of one-mode case with the unit base vector $\xi\in {\mathscr H}$, 
$\|\xi\|=1$.
The $(q,t)$-creation $A_{q,t}^{\dagger} (\xi)$ and 
the $(q,t)$-annihilation $A_{q,t}(\xi)$ operators 
are simply denoted by 
$A_{q,t}^{\dagger}$ and $A_{q,t}$, respectively.
In case of one-mode, the operators ${A_{q,t}}^{\dagger}$ and  $A_{q,t}$
act on the elementary vectors as follows, which can be obtained 
immediately from definitions in Section \ref{sec2}.
In Section \ref{sec:moment}, we discuss combinatorial interpretations 
to these operators.

\begin{lemma}\label{lem:a^*-a}
For $q\in (-1,1)$,  $|q|<t\leq 1$, and $\xi\in\mathscr{H}$ with $\| \xi \|=1$, 
 \begin{align*}
     A_{q,t}^{\dagger} \, & \xi^{\otimes n} = \xi^{\otimes (n+1)}, \; \; n \ge 0,  \quad 
   & A_{q,t}           \, & \xi^{\otimes n} = 
    \begin{cases}
          [n]_{q,t} \, \xi^{\otimes (n-1)}, \; & n \ge 1, \\ 
          0,                                       & n = 0,
    \end{cases} 
\end{align*}
where we adopt $ \xi^{\otimes 0}=\Omega$ as convention.
\end{lemma}
Hence, Lemma \ref{lem:a^*-a} implies following. 
\begin{lemma}\label{lem:qt-n-op}
For $q\in (-1,1)$,  $|q|<t\leq 1$, and $\xi\in\mathscr{H}$ with 
$\| \xi \|=1$, we have 
\begin{equation}\label{lem:number-op}
A_{q,t}^{\dagger}A_{q,t}\xi^{\otimes n}=[n]_{q,t}\xi^{\otimes n},
\ n\geq 1.
\end{equation}
\end{lemma}
It is easy to see that $A_{q,t}^{\dagger}A_{q,t}$ can be identified with the role of 
the so-called number operator $N$ on the Boson Fock space,
and hence it is called $(q,t)$-number operator on ${\mathcal F}_{q,t} ({\mathscr H})$
and denoted by $N_{q,t}:=A_{p,t}^{\dagger}A_{q,t}$ from now on.
In Remark \ref{rem:comb-a^*a}, we give a nice combinatorial interpretation to 
this operator using creation, annihilation, and intermediate cards.

\begin{definition}
For $\lambda > 0$, we consider the bounded self-adjoint operator 
${\bm p}^{\lambda}_{q,t}$
on 
the $(q,t)$-Fock space of one-mode defined by 
\begin{equation}\label{eq:poisson-op}
	{\bm p}^{\lambda}_{q,t}
	:=  N_{q,t} + 
	\sqrt{\lambda}\left(A_{q,t}^{\dagger} + A_{q,t}\right) 
	+  \lambda {\bm 1} ,
\end{equation}
which is our desired model of the $(q,t)$-Poisson type operator 
(random variable) on a noncommutative probability space
$\big( \mathcal{B} \big( {\mathcal F}_{q,t} ({\mathscr H}) \big), 
\varphi \big)$. 
\end{definition}

\begin{theorem}\label{thm:qs-poisson}
For $\xi\in\mathscr{H}$ with $\|\xi \|=1$ and $\lambda>0$, 
the following equality holds: 
\begin{equation}\label{CnSOmega}
C_{n}({\bm p}^{\lambda}_{q,t})\Omega=\sqrt{\lambda^n}\xi^{\otimes n}, \ n\geq 0.
\end{equation}
This equality implies that the probability distribution of ${\bm p}^{\lambda}_{q,t}$ with 
respect to the vacuum expectation is  $\mu_{q,t}^{\lambda}$, 
{\it the $(q,t)$-Poisson distribution} of parameter $\lambda$
\end{theorem}

\begin{proof}
We simply denote ${\bm p}^{\lambda}_{q,t}$ by ${\bm p}$.
We show our claim by induction on $n$. 
It is clear that for $n=0,1$, one can see 
\begin{equation*}
  C_0 ({\bm p}) \, \Omega = {\mathbf 1} \, \Omega = \Omega,           \quad 
  C_1 ({\bm p}) \, \Omega = {\bm p} \, \Omega - \lambda {\mathbf 1} \, \Omega 
                       = \big( \sqrt{\lambda} \,  \xi + \lambda\Omega \big) 
                                        - \lambda {\mathbf 1} \, \Omega
                       = \sqrt{\lambda} \,  \xi.
\end{equation*}
For $n \ge 2$, we assume 
$C_k ({\bm p}) \, \Omega = \sqrt{ \lambda^{k}} \, \xi^{\otimes k}$ 
for $k \le n $. Then it follows that 
 \begin{align*}
   C_{n+1}({\bm p}) \, \Omega 
 &=  \left(
     \left( {\bm p} - (\lambda + [n]_{q,t} ) {\mathbf 1} \right) \, C_n({\bm p}) 
      - \lambda   [n]_{q,t} \, C_{n-1}({\bm p}) 
     \right) \,  \Omega \\
 &=  {\bm p} \, \sqrt{ \lambda^{n}} \, \xi^{\otimes n}
    - (\lambda  + [n]_{q,t} ) \, \sqrt{ \lambda^{n}} \xi^{\otimes n} 
    - \lambda   [n]_{q,t} \, \sqrt{ \lambda^{n-1}} \xi^{\otimes (n-1)} \\ 
 &=  \left( N_{q,t} + \sqrt{\lambda} \, A_{q,t} 
          + \sqrt{\lambda} \, A_{q,t}^\dagger +  \lambda \right) 
                               \sqrt{ \lambda^{n}} \xi^{\otimes n} \\
 &   \qquad \qquad 
    -  \sqrt{ \lambda^{n+2}} \xi^{\otimes n} 
    - [n]_{q,t} \sqrt{ \lambda^{n}  } \xi^{\otimes n} 
    - [n]_{q,t} \sqrt{ \lambda^{n+1}} \xi^{\otimes (n-1)} \\
 &=   [n]_{q,t} \sqrt{ \lambda^{n}} \xi^{\otimes n}
    + [n]_{q,t} \sqrt{ \lambda^{n+1}} \xi^{\otimes n-1}
    + \sqrt{ \lambda^{n+1}} \xi^{\otimes n+1}
    + \sqrt{ \lambda^{n+2}} \xi^{\otimes n} \\
 &   \qquad \qquad 
    - \sqrt{ \lambda^{n+2}} \xi^{\otimes n} 
    - [n]_{q,t} \sqrt{ \lambda^{n}  } \xi^{\otimes n} 
    -  [n]_{q,t} \sqrt{ \lambda^{n+1}} \xi^{\otimes (n-1)} \\
 &=   \sqrt{ \lambda^{n+1}} \xi^{\otimes n+1}.
 \end{align*}
Since $\big\{ C_n ({\bm p}) \big\}_{n \ge 0}$ are self-adjoint operators
with respect to the inner product $\langle \cdot, \cdot \rangle_{q,t}$, we have 
 \begin{align*}
    \big( C_n ({\bm p}) \, C_m ({\bm p}) \, \Omega \, | \, \Omega \big)_{q,t}
& = \big( C_m ({\bm p}) \, \Omega \, | \, C_n ({\bm p}) \, \Omega \big)_{q,t} \\
& = \big( \sqrt{ \lambda^{m}} \xi^{\otimes m} \, | \, 
      \sqrt{ \lambda^{n}} \xi^{\otimes n} \big)_{q,t}\\
& = 0 \; \mbox{ if } \; m \ne n, 
 \end{align*}
which implies 
\begin{equation*}
  \int_{\mathbb{R}}  C_n(t) \, C_m(t) 
                 \, \text{d} \mu_{q,t}^{\lambda}(t) = 0 \; \mbox{ if } \; m \ne n.
\end{equation*}
Therefore, the proof is completed.
\end{proof}

\begin{remark}
(1) One can consider the orthogonal polynomials given in \eqref{eq:qt-op-c} 
as a generalization of the $q$-Charlier polynomials, 
since they include the following well-known examples: 
the $q$-Charlier polynomials of the \textit{Saitoh--Yoshida type} when $t = 1$, 
which appeared in \cite{SY00a, SY00b, An05}, 
and the classical Charlier polynomials \cite{Ch78} in the limit as $t = 1$ and $q \to 1$. 

Moreover, to the best of our knowledge, we have newly found that, for $q = 0$ and $t \in (0,1]$, 
one obtains orthogonal polynomials 
associated with what we call the {\it $t$-free Poisson (Marchenko-Pastur)} distribution 
with parameter $\lambda > 0$, 
which serves as the counterpart 
of the $t$-free Gaussian (semicircular) distribution, 
and reduces to the free Poisson (Marchenko--Pastur) 
distribution \cite{VDN92} when $q = 0$ and $t = 1$. 
Therefore, the polynomials $C_n(x)$ defined in \eqref{eq:qt-op-c} 
can be regarded as a two-parameter deformation of the Charlier polynomials, 
which we refer to as the $(q,t)$-Charlier polynomials.
We note that the term "$t$-free" used here refers to the deformation in \cite{Bl12}, 
and should not be confused with the $t$-deformation examined in \cite{BW01}.

\noindent
(2) On the other hand,  the Poisson type operator 
examined by Asai-Yoshida \cite{AY25} 
is radically different from 
${\bm p}^{\lambda}_{q,t}$ defined above.
Hence, the $(q,t)$-Poisson distribution in this paper 
does not interpolate the Boolean Poisson \cite{SW97} even if $q=0$.
In fact, the Poisson type operator in \cite{AY25} 
is defined by using a different intermediate operator (gauge part) 
from $N_{q,t}$ and in addition 
contains the $s$-deformed operator $k_s$ of identity ${\bm 1}$.
Therefore, the deformed Poisson introduced in this paper 
does not cover deformations of the $s$-free type \cite{AY24} and 
{\it Al-Salam-Carlitz type} \cite{AY25}.

\noindent
(3)
Due to Remark \ref{rem:ej-poisson} and Theorem \ref{thm:qs-poisson},
one can see that ${\bm p}^{\lambda}_{q,t}$ is different from the Poisson type operator 
introduced in \cite{Ej20}. 
\end{remark}

\section{Set Partition Statistics}\label{sec:set-partition}

In our moment formula, the set partitions will be employed as combinatorial objects. 
Here we recall the definition of set partitions and introduce some partition 
statistics for later use.

\begin{definition}
For the set $[n]:= \{1, 2, \ldots , n\}$, a {\it partition of $[n]$} is a collection
$\pi = \{ B_1, B_2, \ldots , B_k \}$ of non-empty disjoint subsets of $[n]$ which 
are called {\it blocks} and whose union is $[n]$. 
For a block $B$, we denote by $|B|$ the size of the block $B$, that is, 
the number of the elements in the block $B$. 
A block $B$ will be called {\it singleton} if $|B| = 1$.
\end{definition}
The set of all partitions of $[n] = \{ 1,2, \ldots , n \}$ will be denoted 
by ${\mathcal P} (\{ 1,2, \ldots , n \})$ or, simply, ${\mathcal P} (n)$.

\subsection{Restricted crossings and nestings:} 
Let $\pi \in {\mathcal P} (n)$ be a partition. 
For elements $e, f \in [n]$, we say that {\it $f$ follows $e$} in $\pi$ 
if $e  <  f$,  $e$ and $f$ belong to the same block of $\pi$, and 
there is no element of this block in the interval $[e, f]$.

\begin{definition}

\begin{itemize}
\item[(1)]
A quadruple $(a, b, c, d )$ of elements in $[n]$ is said to be {\it restricted crossing} of 
$\pi$ if  $c$ follows $a$  in some block of $\pi$ and $d$ follows $b$ in another block 
of $\pi$. 
The statistics $\text{rc}(\pi)$,  {\it the number of restricted crossings of $\pi$}, 
counts the restricted crossings in the partition $\pi$. 
\item[(2)]
A quadruple $(a, b, c, d )$ of elements in $[n]$ is said to be {\it restricted nesting} of 
$\pi$ if  $d$ follows $a$  in some block of $\pi$ and $c$ follows $b$ in another block 
of $\pi$. 
The statistics $\text{rn}(\pi)$,  {\it the number of restricted nestings of $\pi$}, 
counts the restricted nestings in the partition $\pi$. 
\end{itemize}
\end{definition}

\subsection{Graphical Representation:} 
The restricted crossings and nestings have a natural interpretation in the graphic 
line representation of partitions as described, for example, 
in \cite{SU91} and \cite{Bi97}.

Let $\pi$ be a partition in $\mathcal{P}(n)$ and let $B$ be a block of $\pi$. 
If the block $B$ is not singleton (i.e. $|B| \ge 2 $) then we 
write $B = \{ b_1, b_2, \ldots , b_{|B|} \}$. That is, 
$b_{j+1}$ follows $b_{j}$ \; $(j =1, 2, \ldots, |B| - 1)$, and put $b_j$'s on 
$x$-axis.  We will join the points $b_{j}$ and $b_{j+1}$ by an arc above 
the $x$-axis.
Then every restricted crossing appears as {\it a pair of crossing arcs}, and  
also every restricted nesting appears as {\it a pair of nesting arcs}.

\begin{example}\label{ex:4}
(1) $\pi = \big\{ \{1, 3, 4, 7 \}, \{ 2, 5, 10 \}, \{6, 9 \}, \{ 8 \} \big\} 
 \in \mathcal{P}(10)$. 
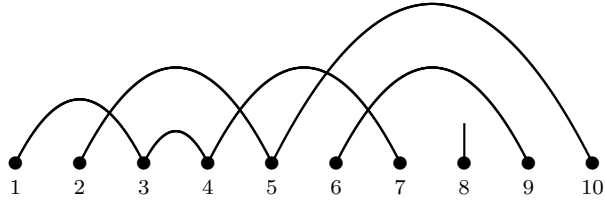
\begin{center}
\setlength{\unitlength}{3pt}
\begin{picture}(72,30)(0, -8)
   \put(  0, 0){\circle*{1.5}} 
   \put(  0,-3){\makebox(0,0){\footnotesize{$1$}}}
   \put(  8, 0){\circle*{1.5}}
   \put(  8,-3){\makebox(0,0){\footnotesize{$2$}}}
   \put( 16, 0){\circle*{1.5}}
   \put( 16,-3){\makebox(0,0){\footnotesize{$3$}}}
   \put( 24, 0){\circle*{1.5}}
   \put( 24,-3){\makebox(0,0){\footnotesize{$4$}}}
   \put( 32, 0){\circle*{1.5}}
   \put( 32,-3){\makebox(0,0){\footnotesize{$5$}}}
   \put( 40, 0){\circle*{1.5}}
   \put( 40,-3){\makebox(0,0){\footnotesize{$6$}}}
   \put( 48, 0){\circle*{1.5}}
   \put( 48,-3){\makebox(0,0){\footnotesize{$7$}}}
   \put( 56, 0){\circle*{1.5}}
   \put( 56,-3){\makebox(0,0){\footnotesize{$8$}}}
   \put( 64, 0){\circle*{1.5}}
   \put( 64,-3){\makebox(0,0){\footnotesize{$9$}}}
   \put( 72, 0){\circle*{1.5}}
   \put( 72,-3){\makebox(0,0){\footnotesize{$10$}}}
  \thicklines
  \drawarc{1}{3}
  \drawarc{3}{4}
  \drawarc{4}{7}
  \drawarc{2}{5}
  \drawarc{5}{10}
  \drawarc{6}{9}
  \drawarc{8}{8}
\end{picture}
\end{center}
Then $\text{rc} (\pi) = 4$ because the partition $\pi$ has four restricted 
crossings, which can be represented by the pairs of crossing arcs  
$\big( [1,3],  [2,5] \big)$, 
$\big( [2,5],  [4,7] \big)$, 
$\big( [4,7],  [5,10] \big)$, and  
$\big( [4,7],  [6,9] \big)$.  
We also find that  $\text{rn} (\pi) = 2$ since there are two 
pairs of nesting arcs  
$\big( [2,5],  [3,4] \big)$ and 
$\big( [5,10],  [6,9] \big)$
as illustrated above.

\noindent
(2) $\pi = \big\{ \{1, 4, 6, 9 \}, \{ 2, 3, 10 \}, \{5 \}, \{ 7, 8 \} \big\} 
 \in \mathcal{P}(10)$. 

\medskip 

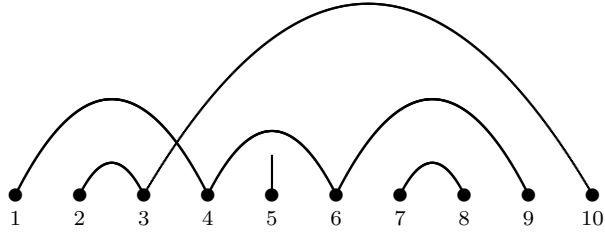
\begin{center}
\setlength{\unitlength}{3pt}
\begin{picture}(72,30)(0, -8)
   \put(  0, 0){\circle*{1.5}} 
   \put(  0,-3){\makebox(0,0){\footnotesize{$1$}}}
   \put(  8, 0){\circle*{1.5}}
   \put(  8,-3){\makebox(0,0){\footnotesize{$2$}}}
   \put( 16, 0){\circle*{1.5}}
   \put( 16,-3){\makebox(0,0){\footnotesize{$3$}}}
   \put( 24, 0){\circle*{1.5}}
   \put( 24,-3){\makebox(0,0){\footnotesize{$4$}}}
   \put( 32, 0){\circle*{1.5}}
   \put( 32,-3){\makebox(0,0){\footnotesize{$5$}}}
   \put( 40, 0){\circle*{1.5}}
   \put( 40,-3){\makebox(0,0){\footnotesize{$6$}}}
   \put( 48, 0){\circle*{1.5}}
   \put( 48,-3){\makebox(0,0){\footnotesize{$7$}}}
   \put( 56, 0){\circle*{1.5}}
   \put( 56,-3){\makebox(0,0){\footnotesize{$8$}}}
   \put( 64, 0){\circle*{1.5}}
   \put( 64,-3){\makebox(0,0){\footnotesize{$9$}}}
   \put( 72, 0){\circle*{1.5}}
   \put( 72,-3){\makebox(0,0){\footnotesize{$10$}}}
  \thicklines
  \drawarc{1}{4}
  \drawarc{2}{3}
  \drawarc{3}{10}
  \drawarc{4}{6}
  \drawarc{5}{5}
  \drawarc{6}{9}
  \drawarc{7}{8}
\end{picture}
\end{center}

\vspace{-15pt}

\noindent
Then $\text{rc} (\pi) = 1$, that is, there is one pair of crossing arcs  
$\big( [1,4],  [3,10] \big)$. While $\text{rn} (\pi) = 5$ since we can find 
the following five pairs of nesting arcs 
$\big( [1,4],  [2,3] \big)$, 
$\big( [3,10],  [4,6] \big)$, 
$\big( [3,10],  [6,9] \big)$, 
$\big( [3,10],  [7,8] \big)$, and  
$\big( [6,9],  [7,8] \big)$.

\end{example}

By using the above set partition statistics, we will derive the combinatorial moment formula of 
the $(q,t)$-deformed Poisson distribution in the next section.

\if10
\subsection{Total Depth of the Blocks by the Last Elements:} 
For our combinatorial formula, we introduce another partition 
statistics related to the last (maximum) elements of the blocks.

For a block $C$ of the partition $\pi \in \mathcal{P}(n)$ we consider 
the first (minimum) element $f_C$ and the last (maximum) element $\ell_C$ in 
the block $C$. In case of singleton it means $f_C = \ell_C$.
For an element $a \in [n]$, we say that the block $C$ {\it covers} $a$ 
if $a$ does not belong to the block $C$ but $a$ is included 
in the interval $[f_C, \ell_C]$.

\begin{definition}
Let $B$ be a block of a partition $\pi$ and $\text{dl}(B)$ count the block 
that covers $\ell_B$ (the last element of $B$), which is called 
{\it the depth of the block $B$ by the last element}.
For a partition $\pi$, the statistics $\text{td}(\pi)$ is defined by 
\begin{equation*}
   \text{td}(\pi) = \sum_{B \in \pi} \text{dl}(B), 
\end{equation*}
which we call {\it the total depth of the blocks by the last elements}.
\end{definition}

\begin{example}
For the partition 
$\pi = \big\{ \{1, 3, 4, 7 \}, \{ 2, 5, 10 \}, \{6, 9 \}, \{ 8 \} \big\} 
\in \mathcal{P}(10)$ in Example 4.3, we put 
$B_1 = \{1, 3, 4, 7 \}$, $B_2 = \{ 2, 5, 10 \}$, $B_3 = \{6, 9 \}$, 
and $B_4 = \{ 8 \}$. 
Then the last element of the block $B_1$ is $7$, which is covered by 
the blocks $B_2$ and $B_3$. Thus $\text{dl}(B_1) = 2$. Similarly we find 
$\text{dl}(B_2) = 0$, $\text{dl}(B_3) = 1$, and  $B_4$ is a singleton and 
covered by $B_2$ and $B_3$, thus $\text{dl}(B_4) = 2$. 
Hence we have $\text{td}(\pi) = 5$.
\end{example}
\fi

\section{Combinatorial moment formula of the $(q,t)$-Poisson distribution}\label{sec:moment}
\if10
Hereafter we will treat the $(q,t)$-Fock space of one-mode with the unit base 
vector $\xi$, and the $(q,t)$-creation $A_{q,t}^{*} (\xi)$ and 
the $(q,t)$-annihilation $A_{q,t}(\xi)$ operators 
are simply denoted by 
$A_{q,t}^{*}$ and $A_{q,t}$, respectively.
We also simply denote the operator ${\bm p}^{\lambda}_{q,t}(\xi)$ by 
$$
    {\bm p}^{\lambda}_{q,t}=
       A_{q,t}^{*} A_{q,t} 
           + \sqrt{\lambda} \big( A_{q,t}^{*} + A_{q,t}  \big)  
              + \lambda \, {\bm 1}.
$$
\fi

Now we are going to investigate the $n$-th moments of the $(q,t)$-Poisson distribution 
of parameter $\lambda$, $\mu^{q,t}_{\lambda}$. 
Namely, we will evaluate the vacuum expectation of the $n$-th power of the 
$(q,t)$-Poisson type operator (random variable) ${\bm p}^{\lambda}_{q,t}$, 
$$
    \Big( \big( {\bm p}^{\lambda}_{q,t}  \big)^n \, \Omega \, \Big| \, \Omega \, \Big)_{q,t}
  =\Big( \big(
        N_{q,t} + \sqrt{\lambda} A_{q,t}^{*} + \sqrt{\lambda} A_{q,t}  
    + \lambda \, {\bm 1} \big)^n \, \Omega \, \Big| \, \Omega \Big)_{q,t}, 
$$
where $N_{q,t} = A_{q,t}^{*} A_{q,t}$ as is mentioned in Lemma \ref{lem:qt-n-op}.

We expand 
$\big( N_{q,t} + \sqrt{\lambda} A_{q,t}^{*} + \sqrt{\lambda} A_{q,t}  
    + \lambda \, {\bm 1} \big)^n$ 
and evaluate the vacuum expectation in a term wise. 
In the expansion, however, we treat all the operators 
$( N_{q,t} ) $, $( \sqrt{\lambda} A_{q,t}^{*} )$, $( \sqrt{\lambda} A_{q,t} )$,
and $( \lambda \, \bm{1} )$  to be noncommutative. 
That is, for example, although $(\lambda \, \bm{1})$ commutes with 
$(\sqrt{\lambda} \, A_{q,t})$ and $(\sqrt{\lambda} \, A_{q,t}^*)$, 
the products 
$(\lambda \, \bm{1}) (\sqrt{\lambda} \, A_{q,t})(\sqrt{\lambda} \, A_{q,t}^*)$
and  
$(\sqrt{\lambda} \, A_{q,t})(\sqrt{\lambda} \, A_{q,t}^*)(\lambda \, \bm{1)}$
should be distinguished each other.  Hence there are $4^n$ terms in the expansion.

We call a term,  a product of the operators 
$( N_{q,t} ) $, $( \sqrt{\lambda} A_{q,t}^{*} )$, $( \sqrt{\lambda} A_{q,t} )$,
and $( \lambda \, \bm{1} )$,  {\it contributor} if it has non-zero 
vacuum expectation.
Associated with a product of $n$ factors 
$$
   Y = Z_n \, Z_{n-1} \cdots Z_2 \, Z_1, 
$$
where 
$Z_k \in \big\{
       ( N_{q,t} ) ,(\sqrt{\lambda} A_{q,t}^{*} ),  ( \sqrt{\lambda} A_{q,t} ), 
       ( \lambda \, \bm{1}) \big\}$ and each factor is numbered from the right, 
we divide $[n] = \{1, 2, \ldots , n \}$ into the following four sets (empty may be allowed):
$$
 \begin{aligned}
     \mathcal{A}_Y & = \big\{\, k \,| \,
                                         z_k = (\sqrt{\lambda} \, A_{q,t} ) \big\}, &
     \mathcal{C}_Y & = \big\{\, k \,| \, 
                                         z_k = (\sqrt{\lambda} \, A_{q,t}^{*} ) \big\}, \\ 
     \mathcal{I}_Y & = \big\{\, k \,| \, 
                                         z_k =  (N_{q,t}) \big\} , & 
     \mathcal{S}_Y & = \big\{\, k \,| \, 
                                         z_k =   ( \lambda \, \bm{1}) \big\}, 
 \end{aligned}
$$
and define the level sequence $\big\{ \ell(k) \big\}_{k = 1}^{n+1}$ by 
$$
   \ell(1) = 0, \quad \; 
   \ell(k + 1) = \ell(k)  + \chi(k). \quad k = 1, 2, \ldots n, 
$$
where $\chi(k)$ is the step function given by 
$$
 \chi(j) = \begin{cases}
                1, \, & \mbox{ if } k \in \mathcal{C}_Y, \\
               -1, \, & \mbox{ if } k \in \mathcal{A}_Y, \\
                0, \, & \mbox{ if } k \in \mathcal{I}_Y \cup \mathcal{S}_Y.
           \end{cases}
$$
Then it can be seen by rather routine argument that if the product $Y$ 
is contributor, that is, 
$\varphi(Y)\ne 0$,
\if10
$\big\langle Y  \, \Omega \, \big| \, \Omega \big\rangle \ne 0$
\fi
then the level sequence 
$\big\{ \ell(k) \big\}_ {k = 1}^{n + 1}$ satisfies the following conditions: 
$$
 \ell(k) \ge 0   \; \, \text{ for } \;  1 \le k \le n, \quad 
 \ell(n+1) = 0, 
 \; \; \text{where if} \; \; k \in \mathcal{I}_Y \; \text{then} \; \ell(k) \ge 1, 
$$
which are equivalent to those for 
$$
  \left( Z_{n} \, Z_{n-1} \cdots \, Z_1 \right) \Omega \in 
                              {\mathbb C} \Omega.
$$
It should be noted that the first two conditions on the level sequence for 
a contributor are known as {\it the Motzkin paths}.

\bigskip

In order to evaluate the vacuum expectation of contributors, we use 
the cards arrangement technique which is similar as in \cite{ER96} for juggling 
patterns. We have already applied this technique for instance in \cite{YY07, Yo20, AY24, AY25},   
but we are now required to assign the different weight of cards to represent the 
restricted nesting and the restricted crossings, appropriately. 
The cards and the weights that we need in this paper are listed below.

\subsection{Creation Cards}\label{subsec:creation-card}
The creation card $C_i$  $(i \ge 0)$  has $i$ inflow lines from the 
left and $(i+1)$ outflow lines to the right, where one new line starts from the middle 
point on the ground level. For each $j \ge 1$, the inflow line of 
the $j$-th level will 
flow out at the $(j+1)$-st level without any crossing. We give the 
weight $\sqrt{\lambda}$ to the card $C_i$. 

The followings are the creation card of the first few levels:
$$
 \begin{array}{ll}
  & \mbox{Level $0$ \Big.} \quad \quad \; \; 
    \mbox{Level $1$ \Big.} \quad \quad \; \; 
    \mbox{Level $2$ \Big.} \ \quad \; \;
    \mbox{Level $3$ \Big.} \\ 
  & \crzer{\sqrt{\lambda}}{C_0}{} \dumcomma \quad 
    \crone{\sqrt{\lambda}}{C_1}{} \dumcomma \quad 
    \crtwo{\sqrt{\lambda}}{C_2}{} \dumcomma \quad 
    \crthr{\sqrt{\lambda}}{C_3}{} \dumcomma \dumcdots
 \end{array} 
$$
The creation card of level $i$ 
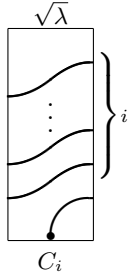
$$
 \begin{array}{l}
    \cicard{\sqrt{\lambda}}{C_i}{} \rsideilong{i} 
 \end{array} 
$$
represents the operation
\begin{equation*}
	\big( \sqrt{\lambda} \, A_{q,t}^{*} \big) \, \xi^{\otimes i} 
     =  \sqrt{\lambda} \, \xi^{\otimes (i+1)}, \quad i \ge 0.
\end{equation*}

\subsection{Annihilation Cards}
The annihilation card $A_i^{(j)}$ $(1 \le j \le i, \, i \ge 1)$ has $i$ 
inflow lines from the left and $(i-1)$ outflow lines to the right. 
On the card $A_i^{(j)}$, only the inflow line of the $j$-th level 
goes down to the middle point on the ground level and ends. 
The lines inflowed at lower than the $j$th level keep their levels. 
The line inflowed at the $\ell( > j)$-th level (higher than the $j$th level) will 
flow out at the $(\ell-1)$-st level (one-decreased level) without any crossing. 
Hence there are $(j-1)$ crossings and $(i - j)$ through lines of non-crossing.
We assign the weight $\sqrt{\lambda} \, t^{i - j} q^{j - 1}$ to 
the card $A_i^{(j)}$, where as is shown below that the parameters 
$q$ and $t$ encode the number of the restricted crossings and the number of the 
restricted nestings, respectively.

The annihilation cards of the first few levels are listed below.
$$
 \begin{array}{ll}
   & \mbox{Level $1$ \Big.} \qquad \qquad \qquad
     \mbox{Level $2$ \Big.} \qquad \qquad \qquad \qquad \quad
     \mbox{Level $3$ \Big.} \\ 
   & \anoneone{\sqrt{\lambda}      }{A_1^{(1)}}{} \dumcomma \qquad 
     \antwoone{\sqrt{\lambda}t     }{A_2^{(1)}}{} \quad 
     \antwotwo{\sqrt{\lambda}   q  }{A_2^{(2)}}{} \dumcomma \qquad  
     \anthrone{\sqrt{\lambda}t^2   }{A_3^{(1)}}{} \quad 
     \anthrtwo{\sqrt{\lambda}t  q  }{A_3^{(2)}}{} \quad 
     \anthrthr{\sqrt{\lambda}   q^2}{A_3^{(3)}}{} \dumcomma \\
   & \qquad \qquad \qquad  \mbox{Level $4$ \Big.} \\ 
   & \anforone{\sqrt{\lambda}t^3   }{A_4^{(1)}}{} \quad 
     \anfortwo{\sqrt{\lambda}t^2q  }{A_4^{(2)}}{} \quad 
     \anforthr{\sqrt{\lambda}t  q^2}{A_4^{(3)}}{}  \quad 
     \anforfor{\sqrt{\lambda}   q^3}{A_4^{(4)}}{} \dumcomma \dumcdots
 \end{array}
$$
The annihilation cards of level $i$,
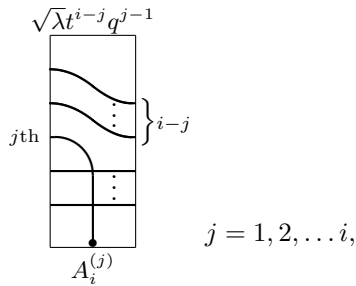
$$
 \begin{array}{rl}
   \aijcard{\sqrt{\lambda} t^{i - j} q^{j - 1}}{A_i^{(j)}}{} \rsideihi{i-j} & \vspace{-40pt} \\
   & j = 1, 2, \ldots i,
 \end{array}
$$
represent the operation 
$$
 \begin{array}{rl}
   \big( \sqrt{\lambda} \, A_{q,t}\big) \, \xi^{\otimes i} 
    & =  \underbrace{\sqrt{\lambda} \, t^{i-1}        \, \xi^{\otimes (i-1)}}_{A_i^{(1)}} + 
            \underbrace{\sqrt{\lambda} \, t^{i-2} q     \, \xi^{\otimes (i-1)}}_{A_i^{(2)}} +
            \underbrace{\sqrt{\lambda} \, t^{i-3} q^2 \, \xi^{\otimes (i-1)}}_{A_i^{(3)}} + \cdots \\ 
    & \qquad \qquad \qquad 
         + \underbrace{\sqrt{\lambda} \,  t q^{i-2} \, \xi^{\otimes (i-1)}}_{A_i^{(i-1)}}  +
            \underbrace{\sqrt{\lambda} \,     q^{i-1} \, \xi^{\otimes (i-1)}}_{A_i^{(i)}}  \\
    & =  \sqrt{\lambda} \, \big( t^{i-1} +  t^{i-2} q + t^{i-3} q^2 
                                          + \cdots + t q^{i-2} + q^{i-1} \big) \, \xi^{\otimes (i-1)} \\
    & =  \sqrt{\lambda} \, [\, i \,]_{q,t} \, \xi^{\otimes (i-1)}, 
         \qquad i \ge 1.
 \end{array}
$$
As one can see, the action $A_{q,t}\xi^{\otimes i}$ can be decomposed 
into $i$ different types of creation cards.

\begin{remark}\label{rem:annihilation-dual}
The total number of 
throughout lines from left to right on $A_i^{(j)}$ is 
$i-1$,  because of $i-1=(i-j)+(j-1)$.
In this sense, the number of restricted nestings
is dual to the number of restricted crossings.
This combinatorial duality on the $(q, t)$-number $[i]_{q,t}$ is explicitly realized by means of a card arrangement, 
which provides a concrete representation of the symmetry between crossings and nestings
of each card.
\end{remark}

\subsection{Intermediate Cards}
The intermediate card $I_i^{(j)}$ $(1 \le j \le i, \, i \ge 1)$ has $i$ 
inflow lines and the same number of outflow lines. 
On the card $I_i^{(j)}$, only the line inflowed at the $j$-th level goes down 
to the middle point on the ground and it will continue as the first 
lowest outflow line. 
The inflow line at the $\ell( < j)$-th level (lower than the $j$-th level) 
will flow out at the $(\ell+1)$-st level (one-increased level), and
the inflow lines of higher than the $j$-th level will keep their levels.
Hence we have $(j-1)$ crossings and there are $(i -j)$ through (horizontal) lines
of non-crossing.
We assign the weight $t^{i-j} q^{j - 1}$ to the card $I_i^{(j)}$, where, 
in the same manner as the annihilation cards, the parameters $q$ and $t$ encode 
the number of the restricted crossings and the number of the restricted 
nestings, respectively.

The intermediate cards of the first few levels are listed below.
$$
 \begin{array}{ll}
   & \mbox{Level $1$ \Big.} \qquad \qquad \quad
     \mbox{Level $2$ \Big.} \qquad \qquad \qquad \qquad \qquad
     \mbox{Level $3$ \Big.} \\ 
   & \nuoneone{1      }{I_1^{(1)}}{} \dumcomma \qquad 
     \nutwoone{t      }{I_2^{(1)}}{} \quad 
     \nutwotwo{    q  }{I_2^{(2)}}{} \dumcomma \qquad  
     \nuthrone{t^2    }{I_3^{(1)}}{} \quad 
     \nuthrtwo{t   q  }{I_3^{(2)}}{} \quad 
     \nuthrthr{    q^2}{I_3^{(3)}}{} \dumcomma \\
   & \qquad \qquad \qquad  \mbox{Level $4$ \Big.} \\ 
   & \nuforone{t^3    }{I_4^{(1)}}{} \quad 
     \nufortwo{t^2 q  }{I_4^{(2)}}{} \quad 
     \nuforthr{t   q^2}{I_4^{(3)}}{}  \quad 
     \nuforfor{    q^3}{I_4^{(4)}}{} \dumcomma \dumcdots
 \end{array}
$$
The intermediate cards of level $i$,
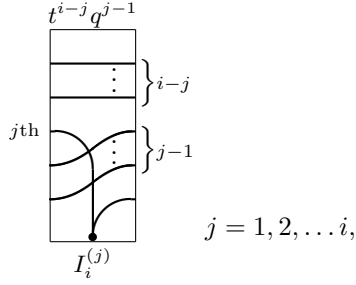
$$
 \begin{array}{rl}
   \nijcard{t^{i-j}q^{j - 1}}{I_i^{(j)}}{} \rsideiji{i-j}{j-1} & \vspace{-40pt} \\
   & j = 1, 2, \ldots i,
 \end{array}
$$
represent the operation 
$$
 \begin{array}{rl}
   \big(N_{q,t} \big) \, \xi^{\otimes i} 
    & = \underbrace{t^{i-1} \, \xi^{\otimes i}}_{I_i^{(1)}} + 
           \underbrace{t^{i-2} q  \, \xi^{\otimes i}}_{I_i^{(2)}} + 
           \underbrace{t^{i-3} q^2  \, \xi^{\otimes i}}_{I_i^{(3)}} + \cdots + 
           \underbrace{t q^{i-2} \, \xi^{\otimes i}}_{I_i^{(i-1)}} + 
           \underbrace{q^{i-1} \, \xi^{\otimes i}}_{I_i^{(i)}} \\
    & =  \big( t^{i-1} +  t^{i-2} q + t^{i-3}q^2 + \cdots 
                      + t q^{i-2} + q^{i-1} \big) \, \xi^{\otimes i} 
      =  [\,i\,]_{q,t} \, \xi^{\otimes i}, 
         \qquad i \ge 1.
 \end{array}
$$
As one can see, 
the action $(N_{q,t}) \, \xi^{\otimes i} $ can be decomposed into $i$ different types of 
intermediate cards.

\begin{remark}
As in the case of annihilation cards, the total number of throughout lines 
from left to right on the card $I_i^{(j)}$ is $i-1$, given by $i-1 = (i-j)+(j-1)$.
In this sense, the number of horizontal non-crossing lines
is dual to the number of crossing lines.  This is a direct consequence of 
Remark \ref{rem:annihilation-dual}.
\end{remark}

\begin{remark}\label{rem:comb-a^*a}
The figure of the lines on the intermediate card $I_i^{(j)}$ can be obtained by the 
composition of those on the annihilation card $A_i^{(j)}$ and the creation card $C_{i-1}$
with identifying the two points on the ground.
$$
 \begin{array}{rl}
   \aijcicard \quad \dumarrow  \qquad \nijcard{}{I_i^{(j)}}{}
 \end{array}
$$
As is mentioned in Lemma \ref{lem:number-op},
this is an interpretation of $N_{q,t}=A_{p,t}^{\dagger}A_{q,t}$
from our card arrangement technical point of view.
On the other hand, the intermediate operator given in \cite{AY25}
does not have such a property.
\end{remark}

\subsection{Singleton Cards}
The singleton card $S_i$  $(i \ge 0)$ has $i$ horizontally parallel lines and the short pole 
at the middle point on the ground. 
We simply assign the weight $\lambda$ to the card $S_i$. 
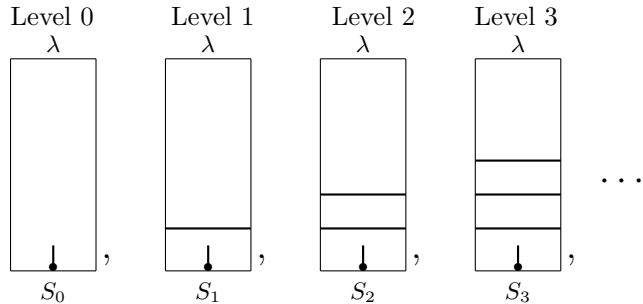
$$
 \begin{array}{ll}
  & \mbox{Level $0$ \Big.} \qquad \; \; 
    \mbox{Level $1$ \Big.} \qquad \; \; 
    \mbox{Level $2$ \Big.} \quad \  \; \;
    \mbox{Level $3$ \Big.} \\ 
  & \kozer{\lambda }{S_0}{} \dumcomma \quad 
    \koone{\lambda }{S_1}{} \dumcomma \quad 
    \kotwo{\lambda }{S_2}{} \dumcomma \quad 
    \kothr{\lambda }{S_3}{} \dumcomma \dumcdots
 \end{array}
$$
The singleton card $S_i$ of level $i$ 
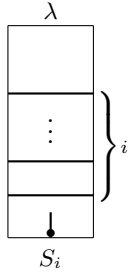
$$
 \begin{array}{l}
    \kicard{\lambda}{S_i}{}\rsideilobig{i}        \vspace{-20pt}
 \end{array}
$$
represents the operation 
$$
   (\lambda \, {\bm{1}}) \, \xi^{\otimes i} = \lambda \, \xi^{\otimes i}.
$$

\subsection{Rules for the Arrangement of the Cards:}
Let $Y = Z_n \, Z_{n-1} \cdots Z_2 \, Z_1$ be a contributor of length $n$, 
and let $\ell(k)$ the level of the $k$-th factor $Z_k$.
Then we arrange the cards along with the following rule:
\begin{itemize}
\item[(1)]{
If $Z_k = (\sqrt{\lambda} \, _{q,t}^*)$, 
we put the creation card of level $\ell(k)$ at the $k$-th site, where 
the only one card $C_{\ell(k)}$ is available.
}
\item[(2)]{
If $Z_k = (\sqrt{\lambda} \, A_{q,t})$, we put the annihilation
 card of level $\ell(k)$ at the $k$-th site, where the cards 
$A_{\ell(k)}^{(j)}$, \;  $1 \le j \le \ell(k)$  are available.
}
\item[(3)]{
If $Z_k = (N_{q,t})$, we put the intermediate card of level $\ell(k)$ at 
the $k$-th site, where the cards 
$I_{\ell(k)}^{(j)}$, \;  $1 \le j \le \ell(k)$  are available.
}
\item[(4)]{
If $Z_k = (\lambda \,{\bm{1}})$, we put the singleton card of level $\ell(k)$
 at the $k$-th site, where the only one card $S_{\ell(k)}$ is available.
}
\end{itemize}

Each card arrangement gives the set partition of $[n]$, 
where the blocks of the partition could be obtained by the concatenation 
of the lines on the cards.  
In this construction, it is easy to find that the creation and the annihilation 
cards correspond to the first (minimum) and the last (maximum) elements in the 
blocks of size greater than equal $2$, respectively, and also that the intermediate cards 
correspond to the intermediate elements in blocks.
Furthermore, {\it the weight of the arrangement} is given by the product of the 
weights of the cards used in the arrangement.

Now we will observe the relation between the weight of the arrangement 
and the set partition statistics:

\bigskip

\noindent
{\bf On the parameter $\bm {\lambda}$:}
\begin{itemize}
\item[$\bullet$]
As is noted in the beginning of Section \ref{sec:moment}, 
the level $\ell(k)$ indicates the $y$-coodinate of the $k$-th step
in a Motzkin path and hence the number of up steps is equal to 
that of downsteps in a Motzkin path.  That is,  we have 
$$
    \#\{ \text{creation cards} \}
  = \#\{ \text{annihilation cards} \}, 
$$
and the parameter $\lambda$ in the product of 
the weights of these cards indicates the number of the blocks of size $\ge 2$.
That is, 
$$
    \big( \sqrt{\lambda} \big)^{\#\{\text{creation cards} \}
                               +\#\{\text{annihilation cards} \} }
 = \lambda^{ \#\{ \text{creation cards}\} }
 = \lambda^{ \# \left\{ 
            B \, \mbox{\footnotesize $|$} \, B \in \pi, \, |B| \ge 2 
            \right\}}.
$$
\item[$\bullet$]
Of course, the parameter $\lambda$ in the product of the weights of singleton 
cards indicates the number of the singletons.
$$
   \lambda^{\#\{ \text{singleton cards} \}}
 = \lambda^{\# \left\{ 
            B \, \mbox{\footnotesize $|$} \, B \in \pi, \, |B| = 1 
            \right\}}.
$$
\end{itemize}
Thus the parameter $\lambda$ in the weight of an arrangement encodes 
the number of blocks of the partition,
$$
    \lambda^{ \# \left\{ 
            B \, \mbox{\footnotesize $|$} \, B \in \pi, \, |B| \ge 2 
            \right\}}
   \lambda^{\# \left\{ 
            B \, \mbox{\footnotesize $|$} \, B \in \pi, \, |B| = 1 
            \right\}}
  = \lambda^{\# \left\{ 
            B \, \mbox{\footnotesize $|$} \, B \in \pi \right\}}
  = \lambda^{|\pi|}.
$$

\bigskip
\noindent
{\bf On the parameters  $\bm{t}$ and $\bm{q}$:}

In a partition obtained by the card arrangement, the blocks are given 
by concatenation of the lines and the elements in the same block are connected 
successively by curves like arcs.

\begin{itemize}
\item[$\bullet$]{
Each ground point on the annihilation cards and the intermediate cards represents 
the right end (larger numbered) site of the arc.  
The non-crossing through lines on the annihilation or the intermediate cards 
are the segments of the arcs, which cover the arc ending at the ground point 
on those cards, because these two arcs will not cross totally.

For instance, in the figure below, the arc $[b,c]$ has the right end (larger numbered) 
site $c$ and  the non-crossing through line on the card at $c$ is the segment 
of the arc $[a,d]$. Then the arc $[a,d]$ covers the arc $[b,c]$ and  
$([a,d], [b,c])$ becomes a pair of nesting arcs.
$$
	\upze \upzc \nestcard \dndz 
$$
That is, each non-crossing through line on the annihilation or the intermediate cards 
corresponds to a restricted nesting in the partition, which is counted by the 
parameter $t$.
Thus the parameter $t$ in the weight of an arrangement encodes the number 
of restricted nestings of the corresponding partition $\pi$, $t^{\text{rn}(\pi)}$.}
\medskip 
\item[$\bullet$]{
Similarly, each crossing through line on the closing or the intermediate cards 
is the segment of the arc, which crosses the arc ending at the ground point 
on those cards, because these arcs will not cross on any other cards.

For instance, in the figure below, the arc $[a,c]$ has the right end (larger numbered) 
site $c$ and the crossing through line on this card is the segment 
of the arc $[b,d]$. Then the arcs $[b,d]$ and $[a,c]$ cross and 
$([a,c], [b,d])$ becomes a pair of crossing arcs.
$$
	\upze \upzc \crosscard \dncz
$$
That is, each crossing through line on the closing or the intermediate cards 
corresponds to a restricted crossing in the partition, which is counted by the 
parameter $q$.
Thus the parameter $q$ in the weight of an arrangement encodes the number 
of restricted crossings of the corresponding partition $\pi$, $q^{\text{rc}(\pi)}$.
}
\end{itemize}

\bigskip 

\begin{example}
The contributor  
$$
 Y = 
  \underbrace{(\sqrt{\lambda} \, A_{q,t})}_{6}
  \underbrace{(\sqrt{\lambda} \, A_{q,t})}_{5}
  \underbrace{(\lambda \,{\bm{ 1}})}_{4}
  \underbrace{( N_{q,t} )}_{3}
  \underbrace{(\sqrt{\lambda} \, A_{q,t}^*)}_{2}
  \underbrace{(\sqrt{\lambda} \, A_{q,t}^*)}_{1}
$$ 
yields the four admissible card arrangements because there are two cards 
available at the sites $3$ and $5$.

We list the admissible card arrangements for $Y$ and their weights 
below. For each arrangement, we also find the corresponding 
partition $\pi \in \mathcal{P}(6)$ and the partition statistics,  
$| \pi |$, $\text{rn}(\pi)$, and $\text{rc}(\pi)$.

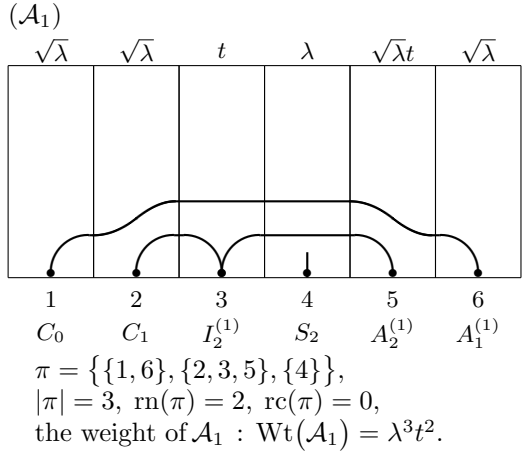
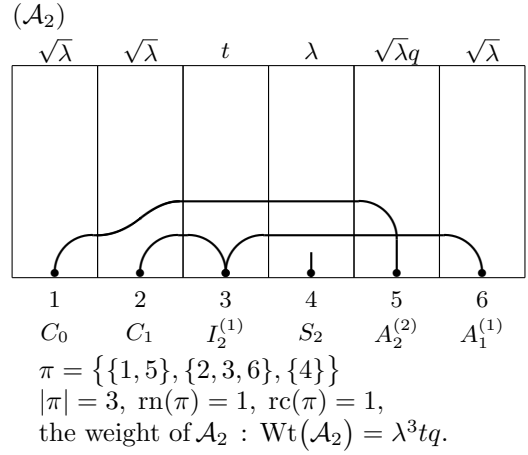
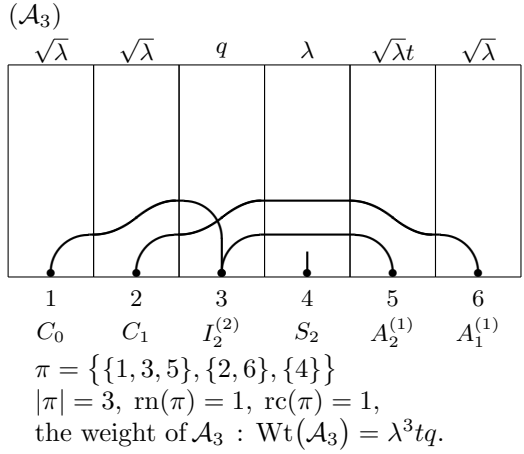
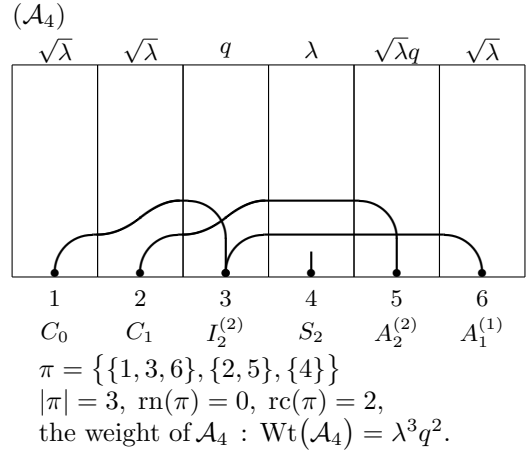
$$
 \begin{array}{l}
  ({\mathcal{A}}_1) \bigg. \\
       \crzer{\sqrt{\lambda}  }{1}{C_0}
       \crone{\sqrt{\lambda}  }{2}{C_1}
    \nutwoone{     t          }{3}{I_2^{(1)}}
       \kotwo{\lambda         }{4}{S_2}
    \antwoone{\sqrt{\lambda} t}{5}{A_2^{(1)}}
    \anoneone{\sqrt{\lambda}  }{6}{A_1^{(1)}} \\
\quad  \pi = \big\{ \{1, 6\}, \{2, 3, 5\}, \{4\} \big\}, \\ 
\quad  |\pi| = 3, \; \text{rn}(\pi) = 2, \; \text{rc}(\pi) = 0, \\
\quad  \mbox{the weight of} \, {\mathcal{A}}_1 \, : \,  
       \mbox{Wt}\big({\mathcal{A}}_1 \big) = \lambda^3 t^2.

 \end{array}
\qquad 
 \begin{array}{l}
  ({\mathcal{A}}_2) \bigg. \\
       \crzer{\sqrt{\lambda}  }{1}{C_0}
       \crone{\sqrt{\lambda}  }{2}{C_1}
    \nutwoone{     t          }{3}{I_2^{(1)}}
       \kotwo{\lambda         }{4}{S_2}
    \antwotwo{\sqrt{\lambda} q}{5}{A_2^{(2)}}
    \anoneone{\sqrt{\lambda}  }{6}{A_1^{(1)}} \\
\quad  \pi = \big\{ \{1, 5\}, \{2, 3, 6\}, \{4\} \big\} \\
\quad  |\pi| = 3, \; \text{rn}(\pi) = 1, \; \text{rc}(\pi) = 1, \\
\quad  \mbox{the weight of} \, {\mathcal{A}}_2 \, : \,  
       \mbox{Wt}\big({\mathcal{A}}_2 \big) = \lambda^3 t q.

 \end{array}
$$

$$
 \begin{array}{l}
  ({\mathcal{A}}_3) \bigg. \\
       \crzer{\sqrt{\lambda} }{1}{C_0}
       \crone{\sqrt{\lambda} }{2}{C_1}
    \nutwotwo{   q           }{3}{I_2^{(2)}}
       \kotwo{\lambda        }{4}{S_2}
    \antwoone{\sqrt{\lambda}t}{5}{A_2^{(1)}}
    \anoneone{\sqrt{\lambda} }{6}{A_1^{(1)}} \\
\quad  \pi = \big\{ \{1, 3, 5\}, \{2, 6\}, \{4\} \big\} \\
\quad  |\pi| = 3,  \; \text{rn}(\pi) = 1, \; \text{rc}(\pi) = 1, \\
\quad  \mbox{the weight of} \, {\mathcal{A}}_3 \, : \,  
       \mbox{Wt}\big({\mathcal{A}}_3 \big) = \lambda^3 t q.
 \end{array}
\qquad 
 \begin{array}{l}
  ({\mathcal{A}}_4) \bigg. \\
       \crzer{\sqrt{\lambda}  }{1}{C_0}
       \crone{\sqrt{\lambda}  }{2}{C_1}
    \nutwotwo{   q            }{3}{I_2^{(2)}}
       \kotwo{\lambda         }{4}{S_2}
    \antwotwo{\sqrt{\lambda} q}{5}{A_2^{(2)}}
    \anoneone{\sqrt{\lambda}  }{6}{A_1^{(1)}} \\
\quad  \pi = \big\{ \{1, 3, 6\}, \{2, 5\}, \{4\} \big\}\\
\quad  |\pi| = 3, \; \text{rn}(\pi) = 0, \; \text{rc}(\pi) = 2, \\
\quad  \mbox{the weight of} \, {\mathcal{A}}_4 \, : \,  
       \mbox{Wt}\big({\mathcal{A}}_4 \big) = \lambda^3 q^2.
 \end{array}
$$
The vacuum expectation of $Y$ is given by 
$$
 \begin{aligned}
	\varphi(Y)
    = \sum_{j=1}^4 \mbox{Wt}\big({\mathcal{A}}_j \big) 
 & = \lambda^3 \,t^2 + \lambda^3 \, t \, q + \lambda^3 \, t  \, q + 
  \lambda^3 \, q^2 \\
 & = \lambda^3 \,  (t + q) (t + q) = \lambda^3 \,  [2]_{q,t} \, [2]_{q,t}. 
 \end{aligned}
$$

\end{example}

\begin{claim}
Every partition in ${\mathcal P}(n)$ can be represented by the admissible card 
arrangement, which is derived from one of those for a certain contributor 
of $n$ factors. 
\end{claim}


\begin{example}
\medskip

\noindent 
(1) \; 
The partition 
$\pi = \big\{ \{1, 3, 4, 7 \}, \{ 2, 5, 10 \}, \{6, 9 \}, \{ 8 \} \big\} 
\in \mathcal{P}(10)$ 
in Example \ref{ex:4} (1) can be obtained by the following admissible card arrangement:
$$
 \begin{array}{l}
       \crzer{\sqrt{\lambda}    }{1 }{C_0}
       \crone{\sqrt{\lambda}    }{2 }{C_1}
    \nutwotwo{ q                }{3 }{I_2^{(2)}}
    \nutwoone{ t                }{4 }{I_2^{(1)}}
    \nutwotwo{ q                }{5 }{I_2^{(2)}}
       \crtwo{\sqrt{\lambda}    }{6 }{C_2}
    \anthrthr{\sqrt{\lambda} q^2}{7 }{A_3^{(3)}}
       \kotwo{\lambda           }{8 }{S_2}
    \antwoone{\sqrt{\lambda} t  }{9 }{A_2^{(1)}}
    \anoneone{\sqrt{\lambda}    }{10}{A_1^{(1)}}
 \end{array}
$$
which is derived from the contributor
$$
 A = 
  \underbrace{(\sqrt{\lambda} \, A_{q,t})  }_{10}
  \underbrace{(\sqrt{\lambda} \, A_{q,t})  }_{9}
  \underbrace{(\lambda \, {\bm 1})             }_{8}
  \underbrace{(\sqrt{\lambda} \, A_{q,t})  }_{7}
  \underbrace{(\sqrt{\lambda} \, A_{q,t}^*)}_{6}
  \underbrace{                    (N_{q,t})}_{5}
  \underbrace{                    (N_{q,t})}_{4}
  \underbrace{                    (N_{q,t})}_{3}
  \underbrace{(\sqrt{\lambda} \, A_{q,t}^*)}_{2}
  \underbrace{(\sqrt{\lambda} \, A_{q,t}^*)}_{1}.
$$ 
The values of partition statistics are given by $|\pi| = 4$, 
$\text{rn}(\pi) = 2$, $\text{rc}(\pi) = 4$ and the weight of card arrangement 
(the product of the weight of the cards) becomes 
$\text{Wt} \big( \mathcal{A}_{\pi} \big) = \lambda^{4} \,t^2 \, q^{4} $.

\bigskip 

\noindent 
(2) \; 
The partition 
$\pi = \big\{ \{1, 4, 6, 9 \}, \{ 2, 3, 10 \}, \{ 5 \}, \{ 7, 8 \} \big\} 
\in \mathcal{P}(10)$ 
in Example \ref{ex:4} (2)  can be obtained by the following admissible card arrangement:
$$
 \begin{array}{l}
       \crzer{\sqrt{\lambda}    }{1 }{C_0}
       \crone{\sqrt{\lambda}    }{2 }{C_1}
    \nutwoone{ t                }{3 }{I_2^{(2)}}
    \nutwotwo{ q                }{4 }{I_2^{(1)}}
       \kotwo{\lambda           }{5 }{S_2}
    \nutwoone{ t                }{6 }{I_2^{(2)}}
    \crtwo{\sqrt{\lambda}       }{7 }{C_2}
    \anthrone{\sqrt{\lambda} t^2}{8 }{A_2^{(1)}}
    \antwoone{\sqrt{\lambda} t  }{9 }{A_2^{(1)}}
    \anoneone{\sqrt{\lambda}    }{10}{A_1^{(1)}}
 \end{array}
$$
which is derived from the contributor
$$
 A = 
  \underbrace{(\sqrt{\lambda} \, A_{q,t})  }_{10}
  \underbrace{(\sqrt{\lambda} \, A_{q,t})  }_{9}
  \underbrace{(\sqrt{\lambda} \, A_{q,t})  }_{8}
  \underbrace{(\sqrt{\lambda} \, A_{q,t}^*)}_{7}
  \underbrace{(N_{q,t})                    }_{6}
  \underbrace{(\lambda \, {\bm 1})               }_{5}
  \underbrace{(N_{q,t})                    }_{4}
  \underbrace{(N_{q,t})                    }_{3}
  \underbrace{(\sqrt{\lambda} \, A_{q,t}^*)}_{2}
  \underbrace{(\sqrt{\lambda} \, A_{q,t}^*)}_{1}.
$$ 
The values of partition statistics are given by $|\pi| = 4$, 
$\text{rn}(\pi) = 5$, $\text{rc}(\pi) = 1$ and the weight of card arrangement 
(the product of the weight of the cards) becomes 
$\text{Wt} \big( \mathcal{A}_{\pi} \big) = \lambda^{4} \,t^{5} \, q$.
\end{example}
\bigskip

Compiling the arguments above, one can obtain the following combinatorial 
moment formula of the $(q,t)$-Poisson distribution:

\bigskip

\begin{theorem}\label{thm:AY-moment}
The $n$-th moment of the $(q,t)$-Poisson distribution
is given by 
\begin{equation*}
	m_{n}(\lambda):=
	\varphi\left( ({\bm p}^{\lambda}_{q,t})^n \right)
    = \! \! 
      \sum_{\begin{subarray}{c} 
               \text{contributor} \, Y \\
               \text{of $n$ factors} 
            \end{subarray}}
	\varphi(Y)
    = \sum_{\pi \in {\mathcal P} (n)} 
       \! \!\lambda^{|\pi|} \, q^{\text{rc}(\pi)} \, t^{\text{rn}(\pi)},
\end{equation*}
where the sum is taking under the set of partitions with $n$ elements, 
$|\pi|$ denotes the number of blocks in $\pi$, 
$\mathit{rc}(\pi)$ is the number of restricted crossings in $\pi$ 
and $\mathit{rn}(\pi)$ is the number of restricted nestings in $\pi$.
\end{theorem}

We list below the first few moments of the $(q,t)$-Poisson distribution :
\begin{align*}
m_1(\lambda)&=1\\
m_{2}(\lambda)&=\lambda^2+\lambda\\
m_{3}(\lambda)&=\lambda^3+(2+t)\lambda^2+\lambda\\
m_{4}(\lambda)&=\lambda^4+(3+t^2+2t)\lambda^3+(3+3t+q)\lambda^2+\lambda.
\end{align*}
Due to the known result of Flajolet \cite{F82},  
the moment generating function of the measure $\mu_{q,t}^{\lambda}$
$$g_{q,t}(z; \lambda):=\sum_{n\geq 0}\Big(\sum_{\pi\in {\mathcal{P}}(n)}q^{rc(\pi)}t^{rn(\pi)}\lambda^{|\pi|}\Big)z^n,$$
has the following continued fraction expression:
\[
g_{q,t}(z; \lambda)=\cfrac{1}{1-\lambda z-\cfrac{\lambda z^2}{1-(\lambda+1)z-\cfrac{\lambda[2]_{q,t}z^2}{1-(\lambda+[2]_{q,t})z-\cfrac{\lambda[3]_{q,t}z^2}{1-(\lambda+[3]_{q,t})z-\cfrac{\lambda[4]_{q,t}z^2}{1-(\lambda+[4]_{q,t})z-\cfrac{\lambda[5]_{q,t}}{{\cdots}}}}}}}
\]

\bigskip
\bigskip
\noindent
{\bf Acknowledgment.}
This work was supported in part by JSPS KAKENHI Grant Numbers JP25K07030 (N. Asai), 
and JP20K03649 and JP25K07046 (H. Yoshida).


\if10
\section*{Conflict of Interest}
On behalf of all authors, the corresponding author states that there is no conflict of interest.

\section*{Data Availability Statement}
This article does not contain any datasets. All results are derived theoretically.
\fi

%

\end{document}